\documentclass[11pt,english]{article}

\usepackage{hyperref}
\hypersetup{
    unicode      = {false},
    pdftoolbar   = {true},
    pdfmenubar   = {true},
    pdffitwindow = {true},
    pdftitle     = {solar},
    pdfauthor    = {Audet et al.},
    pdfsubject   = {solar},
    pdfkeywords  = {solar},
    pdfnewwindow = {true},
    colorlinks   = {true},
    linkcolor    = {blue},
    citecolor    = {blue},
    filecolor    = {black},
    urlcolor     = {blue},
    breaklinks   = {true}
}

\usepackage{amsmath}
\usepackage{amssymb}
\usepackage{multirow}
\usepackage{epsfig}
\usepackage{soul}
\usepackage{colortbl}
\usepackage{curves}
\usepackage{epic}
\usepackage{pgf}
\usepackage{array}
\usepackage{geometry}
\usepackage{times}
\usepackage{xspace}
\usepackage{float}
\usepackage{microtype}
\usepackage{cases}

\def\R{{\mathbb{R}}}
\def\N{{\mathbb{N}}}
\def\Z{{\mathbb{Z}}}
\def\RU{\overline{\R}}
\newcommand{\xx}{\mathbf x}
\newcommand{\nomad}{{\sf NOMAD}\xspace}

\newcommand{\solar}[1]{${\sf solar #1}$\xspace}

\newcommand{\m}{\mathrm{m}} 
\newcommand{\K}{\mathrm{K}} 
\newcommand{\W}{\mathrm{W}} 

\newcommand{\baf}{\mathrm{baf}}
\newcommand{\cold}{\mathrm{cold}}
\newcommand{\cose}{\mathrm{cos}}
\newcommand{\hot}{\mathrm{hot}}
\newcommand{\hf}{\mathrm{hf}}
\newcommand{\pass}{\mathrm{pass}}
\newcommand{\rcv}{\mathrm{rcv}}
\newcommand{\sg}{\mathrm{sg}}
\newcommand{\sh}{\mathrm{sh}}
\newcommand{\spl}{\mathrm{spl}}
\newcommand{\tank}{\mathrm{tank}}
\newcommand{\tub}{\mathrm{tub}}
\newcommand{\trb}{\mathrm{trb}}
\newcommand{\twr}{\mathrm{twr}}
\newcommand{\dry}{\mathrm{dry}}
\newcommand{\wet}{\mathrm{wet}}
\newcommand{\rad}{\mathrm{rad}}
\newcommand{\conv}{\mathrm{conv}}
\newcommand{\bote}{\mathrm{bot}}
\newcommand{\tope}{\mathrm{top}}
\newcommand{\ms}{\mathrm{ms}}

\definecolor{Red}{rgb}{1,0,0}
\definecolor{Green}{rgb}{0,.6,0}
\definecolor{Blue}{rgb}{0,0,1}
\definecolor{Pink}{rgb}{0.7,0,0.7}

\title{ {
\solar{}: A solar thermal power plant simulator for blackbox optimization benchmarking
}
\thanks{This work is partly supported by the NSERC Alliance-Mitacs Accelerate grant  ALLRP~571311-21 (``Optimization of future energy systems'') in collaboration with Hydro-Qu\'ebec.}
}

\author{
\href{mailto:nandres@student.unimelb.edu.au}{Nicolau Andr\'{e}s-Thi\'{o}}\thanks{
  \href{https://www.unimelb.edu.au}{The University of Melbourne},
  \href{mailto:nandres@student.unimelb.edu.au}{\tt nandres@student.unimelb.edu.au}}
\and
    \href{mailto:Charles.Audet@gerad.ca}{Charles Audet}\thanks{
       \href{https://www.gerad.ca}{GERAD} and \href{https://www.polymtl.ca}{Polytechnique Montr\'eal},
          \href{https://www.gerad.ca/Charles.Audet}{\tt www.gerad.ca/Charles.Audet}
  }
  \and
  \href{mailto:diagomartinez.miguel@hydroquebec.com}{Miguel Diago}\thanks{
  \href{https://www.hydroquebec.com/about/}{Hydro-Qu\'ebec}, \href{mailto:diagomartinez.miguel@hydroquebec.com}{\tt diagomartinez.miguel@hydroquebec.com}
  }
  \and
  \href{mailto:aimen.gheribi@polymtl.ca}{{A\"imen~E.} Gheribi} \thanks{
  \href{https://www.crct.polymtl.ca/}{Centre for Research in Computational Thermochemistry}
  and \href{https://www.polymtl.ca}{Polytechnique Montr\'eal}, \href{mailto:aimen.gheribi@polymtl.ca}{\tt aimen.gheribi@polymtl.ca}
  }
  \and
        \href{mailto:Sebastien.Le.Digabel@gerad.ca}{S\'ebastien {Le~Digabel}}\thanks{
        \href{https://www.gerad.ca}{GERAD} and \href{https://www.polymtl.ca}{Polytechnique Montr\'eal},
                  \href{https://www.gerad.ca/Sebastien.Le.Digabel}{\tt www.gerad.ca/Sebastien.Le.Digabel}
}
\and
  \href{mailto:xavier.lebeuf@polymtl.ca}{Xavier Lebeuf}\thanks{
  \href{https://www.gerad.ca}{GERAD} and \href{https://www.polymtl.ca}{Polytechnique Montr\'eal}, \href{mailto:xavier.lebeuf@polymtl.ca}{\tt xavier.lebeuf@polymtl.ca}}
\and
  \href{mailto:mathieu.lemyre-garneau@polymtl.ca}{Mathieu {Lemyre~Garneau}}\thanks{
  \href{https://www.gerad.ca}{GERAD} and \href{https://www.polymtl.ca}{Polytechnique Montr\'eal},
  \href{mathieu.lemyre-garneau@polymtl.ca}{\tt mathieu.lemyre-garneau@polymtl.ca}
  }
  \and
        \href{mailto:christophe.tribes@gerad.ca}{Christophe Tribes}\thanks{
        \href{https://www.gerad.ca}{GERAD} and \href{https://www.polymtl.ca}{Polytechnique Montr\'eal},
                  \href{mailto:christophe.tribes@polymtl.ca}{\tt christophe.tribes@polymtl.ca}
     } 
}

\begin{document}
\maketitle
\thispagestyle{empty}

\thispagestyle{empty}
\noindent
{\bf Abstract:}
{
This work introduces \solar, a collection of  ten optimization problem instances for benchmarking blackbox optimization solvers.
The instances present different design aspects of a concentrated solar power plant simulated by blackbox numerical models.
The type of variables (discrete or continuous), dimensionality, and number and types of constraints (including hidden constraints)
 differ across instances.
Some are deterministic, others are stochastic with possibilities to execute several replications to control stochasticity.
Most instances offer variable fidelity surrogates,
    two are biobjective and one is unconstrained.
The solar plant model takes into account various subsystems: 
	a heliostats field, 
	a central cavity receiver (the receiver), 
	a molten salt  thermal energy storage, 
	a steam generator and 
	an idealized power block.
Several numerical methods are implemented throughout the \solar{} code and most of the executions are time-consuming.
Great care was applied to guarantee reproducibility across platforms.
The \solar{} tool encompasses most of the characteristics that can be found in industrial and real-life blackbox optimization problems,
all in an open-source and stand-alone code.
}



\setcounter{page}{1}

\section{Introduction}

Blackbox optimization (BBO) refers to optimization problems where the objective or constraint functions are not explicitly known or easily computable. 
The term {\em blackbox} refers to the fact that output values can only be obtained via querying the function at the corresponding input points. 
These problems arise in various fields, including engineering design, machine learning, 
    finance and operations research.
An introduction to BBO is found in the textbook~\cite{AuHa2017}, a survey on methodology and software appears in~\cite{CuScVi2017}, a survey on direct-search methods is proposed in~\cite{LaMeWi2019} and hundreds of applications are presented in~\cite{AlAuGhKoLed2020}.

The challenge in BBO lies in the lack of knowledge about the underlying functions. 
Since the internal structure is unknown, traditional optimization methods that rely on explicit mathematical expressions or gradient information cannot be directly applied. 
Instead, algorithmic techniques leverage sampling and exploration to iteratively search the solution space, 
    probing the blackbox for evaluations at different points and using the obtained information to guide the search towards an optimal solution. 

\subsection{Inherent difficulties in BBO}

BBO problems introduce several complexities and considerations. 
First, multiple types of variables might be available: some may be continuous, some discrete and some may even be categorical, i.e., variables which do not satisfy any ordering properties.
In addition, some variables may have an impact on the dimension of the problems.
Terminology for such variables is proposed in~\cite{G-2022-11}.
Second, the evaluation of the functions can be computationally expensive, limiting the number of function evaluations that can be performed. 
Therefore, an optimization algorithm needs to strike a balance between exploration and exploitation to efficiently find good solutions. 
Third, BBO problems often involve noisy or stochastic functions, where the output values may vary even for the same input. 
This necessitates the use of techniques that can handle noise and uncertainty. 
Fourth, in many real applications, 
    the computer simulation used to compute the functions may unexpectedly fail to return valid output.
For example,
 when evaluating a vibration measure of a helicopter rotor blade, 
 approximately 60\% of the simulation calls fail to return a value~\cite{BoDeFrMo0Se98a,BoDeFrSeTo98a}.
These types of constraints are known as hidden constraints~\cite{ChKe00a}.
Some constraints may also return Boolean values, making it difficult for models to approximate them.
Finally,
    there are situations where one wishes to simultaneously optimize more than one objective function.
In this case, one does not wish to find a single solution, 
    but one wishes to obtain a set of non-dominated solutions~\cite{AuSaZg2008a,BiLedSa2020,Custodio2012}.

\subsection{Challenges in BBO benchmarking}

Performance~\cite{DoMo02} and data~\cite{MoWi2009} profiles are now the standard tools for benchmarking derivative-free algorithms.
They allow to agglomerate several optimization runs in simple-to-visualize graphs, and are flexible enough to account for blackbox evaluations. They can also easily be generalized to the constrained and multiobjective cases.
Good benchmarking practices are exposed in~\cite{AuHa2017,Beiranvand2017}.
Accuracy profiles are a related variant described in~\cite{AuHa2017}.
The COCO suite~\cite{HAMTB2016} also provides useful benchmarking tools.
For examples of recent benchmarking studies, see~\cite{PlSah22,RiSa2010}.

Many engineering optimization problems rely on proprietary models and cannot be freely shared with the academic community that develops optimization methods. 
As a consequence, the performance of BBO methods is often assessed on artificial analytic test problems such as the Rosenbrock banana function~\cite{Rose60a}, the Hock and Schittkowski collection~\cite{HoSc1981} or the problems from Mor\'e and Wild~\cite{MoWi2009}.
All of these problems are of crucial importance to nonlinear and derivative-free optimization, but are not representative of difficulties encountered in real BBO applications, in which the objective and constraint functions are not known analytically~\cite{AuHa2017}.
A special issue of {\em Optimization and Engineering}~\cite{AuKok2016} is dedicated to such problems.
The objective of the present work is to introduce a collection of optimization problems as benchmarks for the development of BBO methods. 
The collection is developed from a mechanical engineering perspective and is closer to the problems solved in practice.
Other realistic BBO problems are publicly available. For example,~\cite{FoRe_etal08} compares derivative-free optimization methods  using a set of groundwater supply and hydraulic capture problems. 
A styrene engineering production problem is proposed in~\cite{AuBeLe08}, and~\cite{Abra04,KoAuDe01a} study the optimization of the number and composition of heat intercepts in a load-bearing thermal insulation system.
A pump-and-treat groundwater remediation problem from the Montana Lockwood Solvent Groundwater Plume Site is introduced in~\cite{Matott2011} using the 
{\sf Bluebird} simulator~\cite{bluebird}.

Taken separately, these examples offer desirable characteristics for benchmarking BBO  such as multifidelity, stochasticity, multiobjective and hidden constraints.
In fact, it seems that no work exhibits a realistic application specifically developed for BBO benchmarking.
\solar{} is the first realistic blackbox that gathers all of these characteristics in a single family of instances.
Various optimization of CSP plants are described in~\cite{CoHaNeMa23,LUO2022121716,Manno2020}, but they do not provide benchmarking tools.

\subsection{Contributions}
The main contribution of this work is to provide a family of blackbox benchmarking optimization instances. 
The proposed blackbox gathers most of the desired characteristics of real-world problems, such as 
time-extensive evaluations,
variable fidelity,
stochasticity,
surrogates,
multiobjective,
and several types of constraints, including hidden constraints.

\begin{figure}[htb!]
    \begin{center}
        \includegraphics[width=0.6\textwidth]{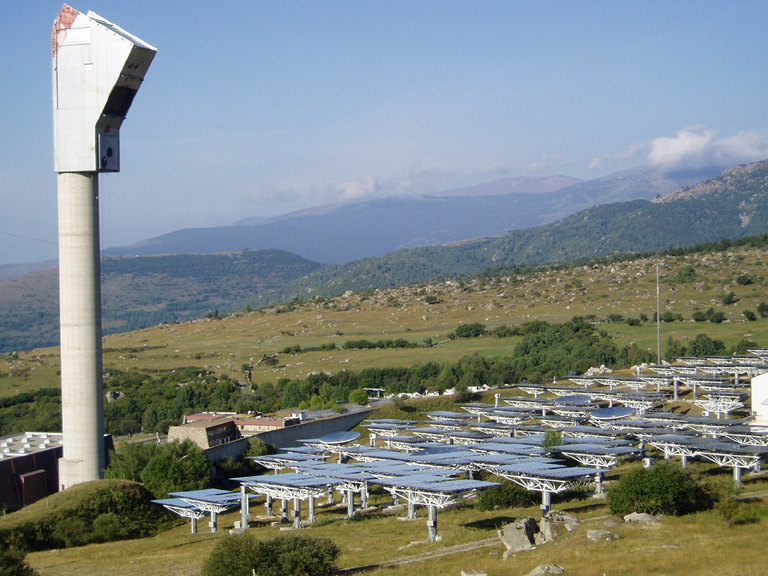} 
    \end{center}
    \caption{ The Th\'emis CSP tower plant in Targasonne, France. \\
    {\scriptsize Source:
    \href{https://commons.wikimedia.org/wiki/File:Themis_2.jpg}
    {\tt https://commons.wikimedia.org/wiki/File:Themis$\_$2.jpg}.}}
    \label{fig-Themis}
\end{figure}

The engineering models proposed  within the blackbox simulates the  operation of a concentrated solar power (CSP)  tower plant with molten salt thermal energy storage.  In this design, a large number of  mirrors, or heliostats, reflects solar radiation on the receiver at the top of a tower.
 The heat collected from the concentrated solar flux is removed from the receiver by a stream of molten salt.
 Hot molten salt is then either used to feed thermal power to a conventional power block,  or stored in  an insulated tank for deferred use.  The Th\'emis CSP power plant~\cite{drouot1984themis} (see Figure~\ref{fig-Themis}), in France, was the first built with this design.  Gemasolar, in Spain, was the first CSP tower plant to produce electrical power on a ``24/7'' basis using molten salt thermal energy storage~\cite{burgaleta2011gemasolar}.


The source code of the models  is freely available under the LGPL license at \href{https://github.com/bbopt/solar}{\url{https://github.com/bbopt/solar}}, along with some starting points, best known values, and the logs provided by some solvers such as \nomad~\cite{nomad4paper,Le09b}.
\solar{} is written in standard {\sf C++}, and compiles on most platforms and is guaranteed to reproduce the exact same outputs independently of the platform.
It also provides variable-fidelity static surrogates,
deterministic or stochastic outputs,
and the control over its number of replications, leading to reduced variance in the output at the cost of more time-expansive simulations.
The code is parametrized so that a total of 10 BBO problems instances with various complexity levels is provided.
These instances are  single- or bi-objective,  involve between 5 and 29 integer, real or categorical variables, and implement between 5 and 17 binary, integer or continuous constraints.
In addition, hidden constraints are present as some evaluations may fail to compute.


\subsection{Organization}

The paper is structured as follows: Section~\ref{sec-plant} provides a short literature review on solar plants, and a high-level description of the main components of the power plant model.
The complete model description is available  separately, in the MSc thesis~\cite{MScMLG}.
Section~\ref{sec-instances} describes the collection of test problems
from the point of view of BBO benchmarking.
Section~\ref{sec-thepackage} describes the \solar{} package, its many features, and provides illustrative examples of typical benchmarking situations.
Concluding remarks follow.

\section{Description of the solar thermal power plant simulator}
\label{sec-plant}

This section gives a short literature review on CSP plants and introduces the technologies that are simulated when different instances are queried, as well as the mathematical/numerical models used to represent them.

\subsection{A short literature review on solar plants}
\label{sec-litrev}


Solar energy has emerged as a pivotal solution to address challenges associated with fossil fuel depletion, environmental degradation, and climate change. Thermal energy, derived from solar radiation, represents a potentially abundant energy source. Harnessing and storing energy from solar radiation primarily involve two methods: photovoltaic and thermal energy storage.

Photovoltaics convert solar radiation into electrical energy, which can be used immediately or stored using batteries. 
When integrated into the grid, inverters are employed to convert direct current into alternating current. 
However, the intermittent and variable nature of solar radiation introduces substantial unpredictability into its utilization. 
Therefore, the development of efficient energy storage systems is imperative to fully realize its benefits. 
Thermal energy storage technologies entail the transfer of energy to a material, leading to an increase in its total enthalpy.
This stored energy can be deployed as needed to meet grid power demands. Ultimately, the stored energy is transferred to a heat-transfer fluid, often water, using either sensible heat (heat transfer without a change in state) or latent heat interactions. 

One example of sensible heat storage technology is the Solar Electric Station IX in California's Mojave Desert, which utilizes a Phase Change Material (PCM). 
PCMs enable energy storage or release through a first-order phase transition, and the enthalpies of these transitions, known as latent heat, quantify the heat interactions. 
Four possible phase transformations can be used for heat storage: solid-solid, solid-liquid, solid-gas, and liquid-gas reactions. 
PCMs that undergo solid-liquid phase transitions offer the highest energy density variation, providing gigajoules of energy compared to kilojoules for liquid-gas transitions. 
Wei et al.~\cite{wei2018} have reviewed the selection principles for PCMs, considering material properties like heat of transition, melting mechanism, thermal conductivity, heat capacity, volume change upon phase change, vapor pressure, reactivity, operating temperature constraints, and heat transfer equipment design. 

Thermophysical properties are essential in PCM design. Numerous materials can serve as PCMs, covering a wide range~\cite{farid2004,kenisarin2007,kenisarin2010,sharma2009}. They are often categorized based on their chemical properties, with classifications for low-temperature thermal energy storage by Abhat~\cite{abhat83} and for medium and high-temperature PCMs in~\cite{kenisarin2010}. 
In~\cite{GHERIBI2020110506}, a comprehensive analysis of the intricate interplay between thermophysical, thermal, economic, and corrosion properties identifies optimal PCMs for sustainable CSP plants. The study identified criteria including exceptional heat storage capacity, optimal thermal conductivity in solid and liquid phases, and high heat capacity per unit volume. 

\subsection{Introduction to the simulator}

Several CSP subsystems are involved in the transformation of direct solar irradiation to electrical power.
Table~\ref{tab-nomenc}
details the quantities and notations used through this work.
Note that these quantities do not necessarily correspond to the {\em optimization} variables described in Section~\ref{sec-instances}.
    
\begin{table}[ht!]
    \renewcommand{\tabcolsep}{3pt}
    \centering
    \begin{tabular}[t]{lll}
        \textbf{Name}  & \textbf{Quantity} & \textbf{Units} \\
        \hline
        \hline
$D$ & outer diameter & $\m$ \\ 
$d$ & inner diameter & $\m$ \\ 
$H$ & height & m \\ 
$L$ & length & m \\ 
$P$ & parasitic load & k$\W$ \\ 
$Q$ & heat transfer & k$\W$ \\ 
$R$ & radius & $\m$ \\ 
$ST$ & type of steam turbine & $-$ \\ 
$t$ & thickness & $\m$ \\ 
$W$ & width & $\m$ \\ 
$\eta$ & efficiency & $-$  \\ 
$\theta$ & angle & deg \\ 
$\tau$ & transmissivity & $-$  \\ 
        \hline
    \end{tabular}
%
\quad
%
    \centering
    \begin{tabular}[t]{lll}
        \textbf{Subscript}  & \textbf{Definition} \\
        \hline
        \hline
$\baf$ & baffles \\
$\bote$ & bottom \\
$\cold$ & cold tank \\
$\conv$ & convection \\
$\cose$ & cosine effect \\
$\hot$ & hot tank \\
$\hf$ & heliostats field \\
$\ms$ & molten salt \\
$\pass$ & passes in steam generator \\
$\rad$ & radiation \\
$\rcv$ & receiver \\
$\sg$ & steam generator \\
$\sh$ & shell \\
$\spl$ & spillage \\
$\tank$ & tank ($\hot$ or $\cold$) \\
$\tub$ & tubes \\
$\trb$ & turbine \\
$\twr$ & tower \\
        \hline
    \end{tabular}
\caption{Notations for various quantities and signification of subscripts.
Subscripts ``$\dry$'', ``$\tope$'' and ``$\wet$'' are obvious.}
\label{tab-nomenc}
\end{table}

\subsection{System dynamics}

\begin{figure}
    \includegraphics[width=\linewidth]{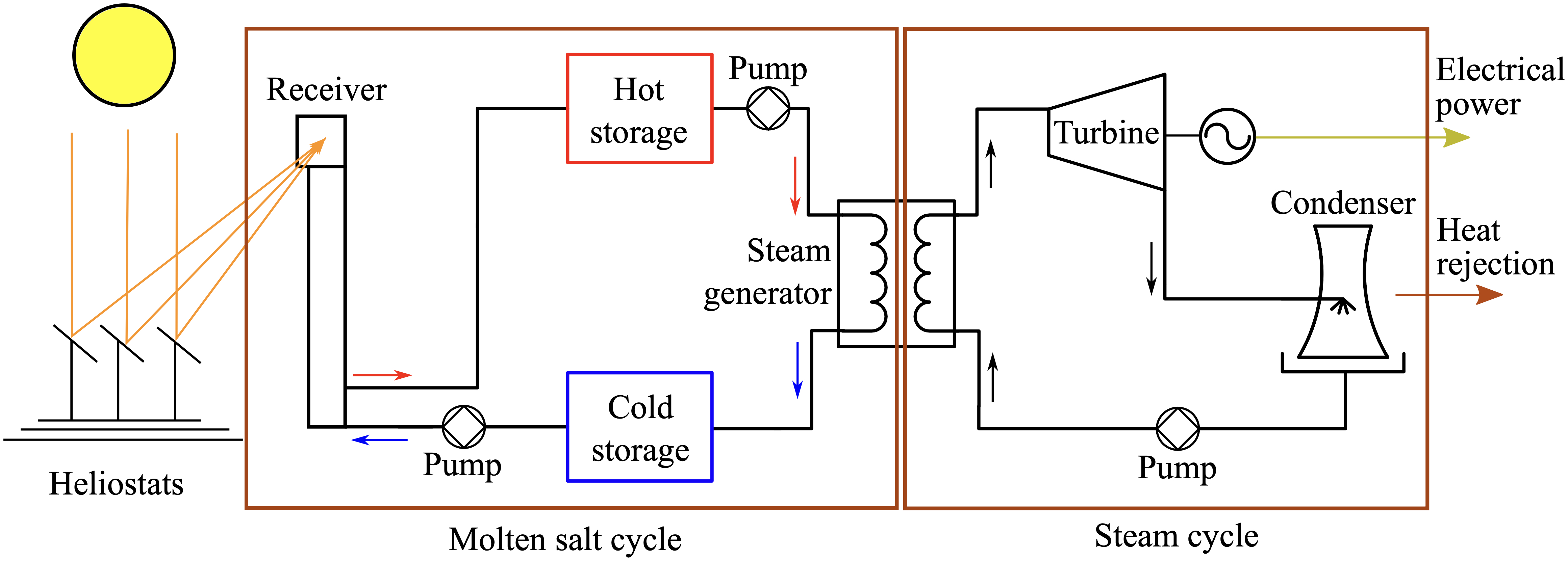}
    \caption{Cycle of the concentrated solar power (CSP) plant.}
    \label{fig-PhysDynamic}
\end{figure}

 The plant is designed to transform heat recovered from concentrated sunlight into electrical power. Due to the Second Law of Thermodynamics and to non-idealities in each subsystem, part of the input heat is rejected during the process to the environment.
Figure~\ref{fig-PhysDynamic}  provides an overview of the complete cycle. Losses due to non-idealities are not represented but are accounted for in all components except the steam generator

 The heliostats field concentrates sunlight on the receiver, at the top of the tower. As the receiver heats up, thermal power is extracted by raising the temperature of a {flow} 
 of molten salt pumped through its structure.
 The resulting hot molten salt is then directed to  a hot storage  tank until it is needed to drive the power block. 

 When electrical power is produced, hot molten salt is pumped though the steam generator to operate the power block. Heat is transferred to a {flow} 
 of water on the other side of the steam generator which is transformed to superheated steam.
Cold molten salt is  recovered in the cold storage  tank  where it remains until it is pumped through the receiver again.

 The power block is defined as a simple Rankine cycle with a single steam turbine.
 Superheated high-pressure steam drives a turbine coupled to an electrical generator. Low-pressure steam is then condensed and pumped back as liquid water though the steam generator.
  
\subsection{Heliostats field}

The optical  subsystem consists of an array of heliostats  (i.e., sun tracking mirrors)  that reflect  sunlight on the receiver. The design parameters used to define the heliostats field is given in Table~\ref{tab-HeliostatsField}.
 All of these parameters, except the latitude, are considered as optimization variables in Section~\ref{sec-vars}.

\begin{table}[ht!]
\centering
\begin{tabular}{lllr}
\hline
\multirow{2}{*}{Symbol}	& \multirow{2}{*}{Quantity}& \multirow{2}{*}{Unit}& {Optim. var. in}\\
&&& {Tables~\ref{tab-vars1} and~\ref{tab-vars2}}\\
\hline
\hline
$\phi$ & latitude & $\deg$ & $-$ \\
$ W_{\hf}$  & heliostats width  & $\m$ & $x_2$	\\
$ L_{\hf}$ & heliostats length	& $\m$	& $x_1$	\\
$ H_{\twr}$ & height of tower	& $\m$	& $x_3$	\\
$ N_{\hf}$ & number of heliostats to fit in the field & $-$	& $x_6$ \\
$\theta_{\hf}$ & angular width of the heliostats field & \multirow{2}{*}{$\deg$} & \multirow{2}{*}{$x_7$} \\
			   & on each side of the N-S axis & & \\
$R_{\hf}^{\min}$ & minimum distance between heliostats & \multirow{2}{*}{proportion of $H_{\twr}$} & \multirow{2}{*}{$x_8$} \\
			     & and tower					       &   & 	 	\\
$R_{\hf}^{\max}$ & maximum distance between heliostats	& \multirow{2}{*}{proportion of $H_{\twr}$} & \multirow{2}{*}{$x_9$}\\
			& and tower						&	& \\
$H_{\rcv}$ & receiver aperture height & $\m$ & $x_1$, $x_4$	\\
$W_{\rcv}$ & receiver aperture width	& $\m$	& $x_2$, $x_5$	\\
\hline
\end{tabular}
\caption{List of design parameters for the heliostats field. The nomenclature of Table~\ref{tab-nomenc} is used.
The latitude $\phi$ is a parameter fixed during optimization.}
\label{tab-HeliostatsField}
\end{table}

\subsubsection{Generating the heliostats field}

\begin{figure}[ht!]
    \centering
    \includegraphics[width=0.49\textwidth]{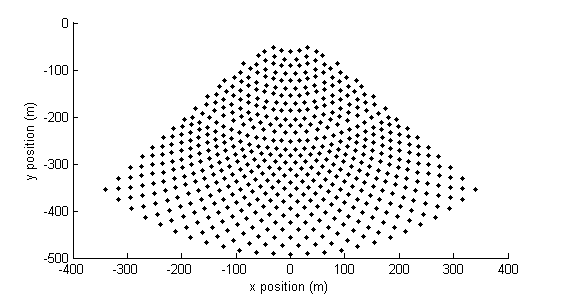} 
    \includegraphics[width=0.49\textwidth]{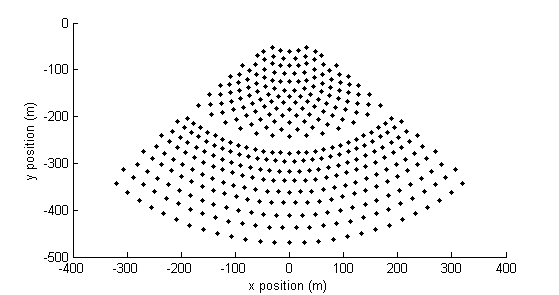}
    \caption{Examples of heliostats field layouts  for tower heights of 120~$\m$  (left) and 70~$\m$  (right).  The tower is located at (0,0).}
    \label{fig-hfLayout}
\end{figure}

 The heliostats are laid on a radially staggered grid that prevents blocking losses between them~\cite{Siala2001}. The grid is calculated as a function of individual heliostat dimensions ($W_{\hf}$ and $L_{\hf}$) and tower height ($H_{\twr}$).
Figure~\ref{fig-hfLayout} shows two examples of radially staggered grid layouts for identical heliostat dimensions but different tower heights. 

\subsubsection{Heliostats layout}

Once  the grid layout  is determined, each position is rated according to  the average optical efficiency $\eta$  that a potential heliostat would have on that location  during the time window simulated.  The optical efficiency of a heliostat is  considered as the product of its cosine efficiency ($\eta_{\cose}$), the total atmospheric transmissivity from  its position (constant $\tau$), and  its spillage efficiency ( $\eta_{\spl}$). Both $\eta_{\cose}$ and $\eta_{\spl}$ depend on the orientation of each heliostat, which is itself dependent on the sun position.  Both losses are approximated in~\cite{MScMLG} as a function of an analytical approximation of the Sun position~\cite{duffie2013solar}. Shadowing effects are not considered during this step, but are be accounted for when calculating the overall field's performance. 
The  instantaneous efficiency  of a potential heliostat position is defined as:
\begin{equation*}
     \eta = \eta_{\cose} \eta_{\spl} \tau
\end{equation*}

\begin{figure}[ht!]
\centering
\includegraphics[width=0.49\textwidth]{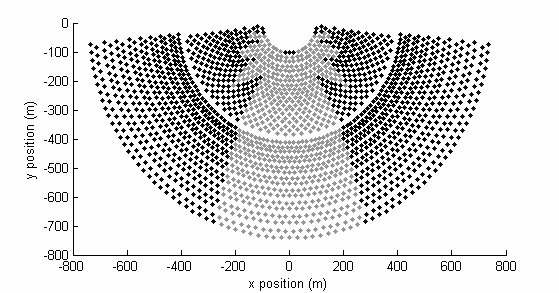} 
\includegraphics[width=0.49\textwidth]{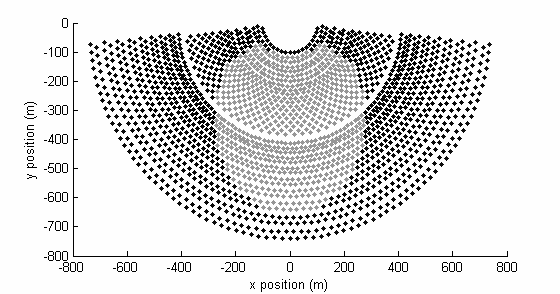}
\caption{Examples of positions selected (in light gray) for receiver aperture widths of 3~$\m$ (left) and 15~$\m$ (right).}
\label{fig-hfSelection}
\end{figure}

The actual heliostats field is generated by occupying the first $N_{\hf}$ grid positions with the highest average optical efficiency for the given receiver aperture and tower height.
 For illustration, Figure~\ref{fig-hfSelection} shows how the arrangement of 700 heliostats on the same spatial grid of 1,960 points varies with the receiver aperture width ($W_{\rcv}$).
 As the aperture narrows, the algorithm selects heliostats closer to the North-South axis to minimize spillage.
 For wider apertures, the selection is dictated by cosine efficiency and atmospheric attenuation.

\subsubsection{Sun radiation model}





 The power delivered to the receiver is estimated with a ray-tracing procedure corrected for spillage, attenuation and reflectivity losses. The sun is approximated as a distant point source  so that rays are parallel to one another when they reach the field. Direct solar irradiance at the Earth's surface is assumed to be 1~k$\W/\m^2$.  Attenuation losses are calculated for a clear day~\cite{ballestrin2012solar}.  Heliostats are modelled as flat surfaces with no aiming or tracking errors.  Shading and cosine losses are inherently accounted for during ray-tracing.
 As a result of these simplifications, the predictions delivered by the model may overestimate the actual performance of the CSP plant.

%
%

\subsection{The receiver}

 The receiver is modelled as a cavity with a rectangular aperture of sides $H_{\rcv}$ and $W_{\rcv}$. Molten salt is heated as it flows through an array of tubes laid on the inner surface of the cavity which are exposed to concentrated solar irradiation.
Table~\ref{tab-CentralReceiver}  lists the design parameters  used to model the receiver.

\begin{table}[ht!]
\centering
\begin{tabular}{llll}
\multirow{2}{*}{\bf Symbol} & \multirow{2}{*}{\bf Quantity} & \multirow{2}{*}{\bf Unit} & {\bf Optim. variables in} \\
		  	                &					            &		                    & {\bf Tables~\ref{tab-vars1} and~\ref{tab-vars2}} \\
\hline
\hline
 $H_{\rcv}$	                & aperture height	            & $\m$                      & $x_1$, $x_4$\\
 $W_{\rcv}$ 	            & aperture width	            & $\m$                      & $x_2$, $x_5$\\
 $N_{\rcv,\tub}$	        & number of tubes 	     	    & $-$                       & $x_4$, $x_7$, $x_{10}$, $x_{11}$, $x_{16}$ \\
 $D_{\rcv}$ 		        & tubes outer diameter	        & $\m$                      & $x_7$, $x_{10}$, $x_{13}$, $x_{14}$, $x_{19}$\\
 $d_{\rcv}$  		        & tubes inner diameter	        & $\m$                      & $x_6$, $x_9$, $x_{12}$, $x_{13}$, $x_{18}$\\
 $t_{\rcv}$ 		        & thickness of insulation    	& $\m$                      & $x_5$, $x_8$, $x_{11}$, $x_{12}$, $x_{17}$ \\
\hline
\end{tabular}
\caption{List of design parameters for the receiver unit. The nomenclature of Table~\ref{tab-nomenc} is used.
$H_{\rcv}$ and $W_{\rcv}$ also appear in Table~\ref{tab-HeliostatsField} since they impact the design of the heliostats field.
}
\label{tab-CentralReceiver}
\end{table}

 Heat lost back to the environment by reflection, re-radiation and conduction through the receiver wall  is modelled with an iterative procedure described by Li~et~al.~\cite{li2010thermal}. The absorbed heat rate  is determined  as the difference between the receiver thermal losses and the solar flux concentrated by the heliostats field on its aperture.

\subsection{Thermal energy storage}

The two thermal  energy storage units, namely the cold tank and the hot tank, are large  cylindrical stainless steel tanks  where molten salt is kept until  pumped to the power block or to the solar receiver.
Both tanks have the same diameter and the cold storage height is 1.2 times larger than
the hot storage height.
Both tanks are protected by a layer of  thermal insulation.
Table~\ref{tab-ThermalStorage} shows the list of design parameters that define the thermal storage  component.

\begin{table}[ht!]
\centering
\begin{tabular}{llll}
\multirow{2}{*}{\bf Symbol} & \multirow{2}{*}{\bf Quantity} & \multirow{2}{*}{\bf Unit}     & {\bf Optim. variables in}     \\
		  	                &				                &                               & {\bf Table~\ref{tab-vars1}} \\
\hline
 $t_{\tank}$                & tanks insulation thickness                    & $\m$          & $x_4$, $x_5$, $x_{13}$, $x_{14}$\\
 $H_{\tank}$                & height of the interior of the storage tanks   & $\m$          & $x_2$, $x_{11}$\\
 $d_{\tank}$                & diameter of the interior of the storage tanks & $\m$          & $x_3$, $x_{12}$\\
\hline
\end{tabular}
\caption{List of design parameters for the thermal storage units.
The nomenclature of Table~\ref{tab-nomenc} is used.
The subscript ``$\tank$'' can either be ``$\hot$'' or ``$\cold$'' and denotes whether the term refers to the hot or the cold tank e.g. $t_{cold}$ denotes the cold tank insulation thickness.
The geometry of tanks is constrained by $d_{\cold}=d_{\hot}$ and $H_{\cold}=H_{\hot}\times 1.2$.
}
\label{tab-ThermalStorage}
\end{table}

\begin{figure}[ht!]
\setlength{\unitlength}{0.8cm}
\includegraphics[width=8cm, height=8cm]{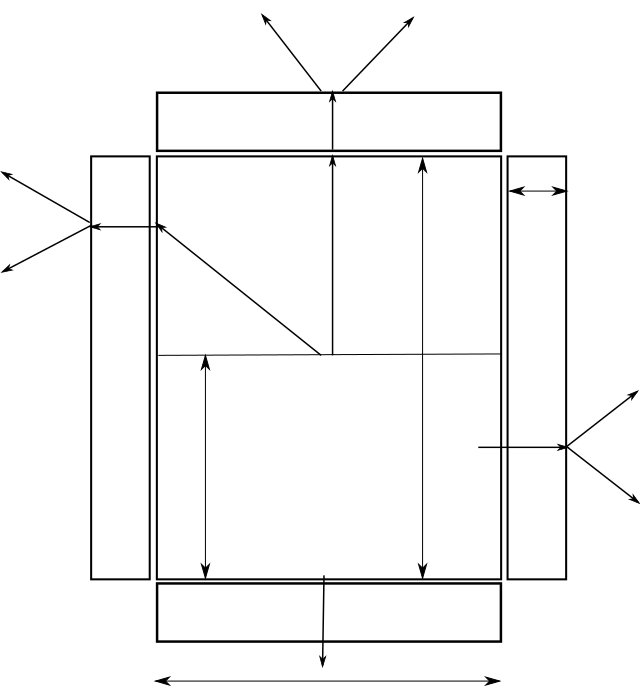}
\centering
\begin{picture}(0,0)(10,0)
\put(4,10){\makebox(0,0){$q_{\tope,\rad}$}}
\put(6,10){\makebox(0,0){$q_{\tope,\conv}$}}
\put(0,7.8){\makebox(0,0){$q_{\dry,\rad}$}}
\put(0,5.8){\makebox(0,0){$q_{\dry,\conv}$}}
\put(10,4.75){\makebox(0,0){$q_{\wet,\rad}$}}
\put(10,2.5){\makebox(0,0){$q_{\wet,\conv}$}}
\put(9.5,7.5){\makebox(0,0){$t_\tank$}}
\put(5.8,4){\makebox(0,0){$H_\tank$}}
\put(5,-0.1){\makebox(0,0){$d_\tank$}}
\put(3.75,3.5){\makebox(0,0){$L_{\ms}$}}
\put(6,6){\makebox(0,0){$q_{\ms \to \tope}$}}
\put(4,6.3){\makebox(0,0){$q_{\ms \to \dry}$}}
\put(5,2){\makebox(0,0){$q_{\bote}$}}
\end{picture}
\caption{Thermal storage model parameters definition and loss processes.
The nomenclature of Table~\ref{tab-nomenc} is used.
The subscript ``$\tank$'' can either be ``$\hot$'' or ``$\cold$''.
}
\label{fig-SimpleHeatStorage}
\end{figure}
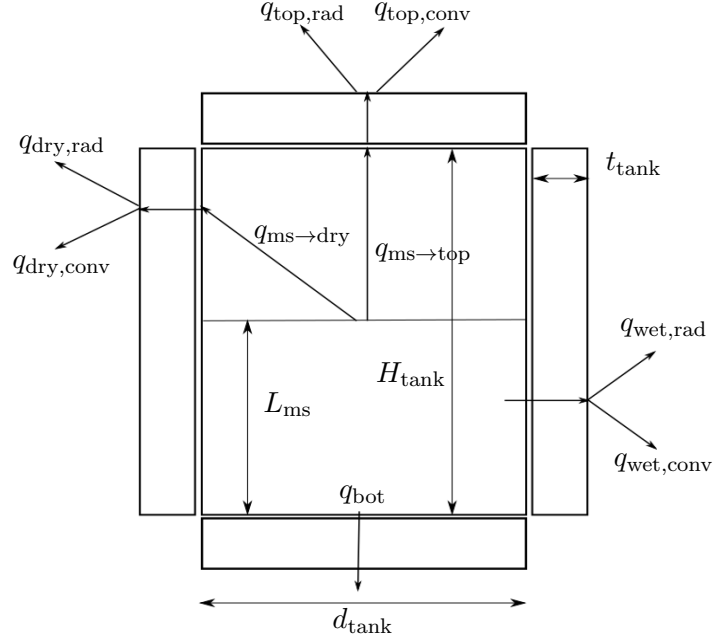

 Thermal losses during storage are  estimated as a function of the  molten salt level and temperature  with the heat loss model proposed by Zaversky~et~al.~\cite{zaversky2013transient}. Figure~\ref{fig-SimpleHeatStorage} illustrates the heat transfer modes considered  and defines the  tanks design parameters.
 An identical insulation thickness for the wall and the ceiling  surfaces is assumed.  Natural convection and radiation patterns within the air and molten salt volumes are neglected.  External radiation and forced convection  losses to the atmosphere  are considered assuming a constant wind  speed of 6~$\mathrm{\m.s^{-1}}$  and normal atmospheric conditions.
 There are three channels of thermal losses from the molten salt when the tank is not at full capacity: conduction through the tank’s floor, conduction through the wet part of the cylindrical wall, and radiation from its top surface to the ceiling and the dry part of the cylindrical wall.

\subsection{Heat exchanger (steam generator)}
\label{sec-heat-exchanger}

 In general, commercial CSP plants obtain superheated steam from the combined action of a pre-heater, boiler and super-heater (e.g.,~\cite{shams2017}). None of these subsystems are explicitly modelled here. Instead, the steam generation unit is modelled only as a shell-and-tubes molten salt to water heat exchanger using the Effectiveness-NTU method~\cite{incroperafundamentals},
 with the use of baffles that support the tubes, and
 the baffles cut corresponds to the portion of 
the shell inside diameter that is not covered by the baffle.
  Table~\ref{tab-SteamGen} lists the design parameters associated with the heat exchanger.

\begin{table}[ht!]
\begin{center}
\begin{tabular}{llll}
\multirow{2}{*}{\bf Symbol} & \multirow{2}{*}{\bf Quantity}   & \multirow{2}{*}{\bf Unit} & {\bf Optim. variables}            \\
                            &                                 &                           & {\bf in Table~\ref{tab-vars2}} \\
\hline
 $L_{\sg}$                  & length of tube passes           & $\m$    & $x_{12}$, $x_{21}$ \\
 $d_{\sg}$                  & inner diameter of tubes         & $\m$    & $x_{13}$, $x_{22}$\\
 $D_{\sg}$                  & outer diameter of tubes         & $\m$    & $x_{14}$, $x_{23}$\\
 $H_{\sg,\baf}$             & baffles cut       & ratio of shell width  & $x_{15}$, $x_{24}$\\
 $N_{\sg,\baf}$             & number of baffles               & $-$     & $x_{16}$, $x_{25}$\\
 $N_{\sg,\tub}$             & number of tubes                 & $-$     & $x_{17}$, $x_{26}$\\
 $N_{\sg,\sh,\pass}$        & number of shell passes          & $-$     & $x_{18}$, $x_{27}$\\
 $N_{\sg,\tub,\pass}$       & number of tube passes per shell & $-$     & $x_{19}$, $x_{28}$\\
\hline
\end{tabular}
\caption{List of design parameters for the shell-and-tubes heat exchanger model.
The nomenclature of Table~\ref{tab-nomenc} is used.}
\label{tab-SteamGen}
\end{center}
\end{table}

The heat exchanger model has two main functions. The first is to determine the flow of hot molten salt necessary to produce the required amount of pressurized steam for the turbine. For any simulation interval, the computed molten salt flow is dependent on the heat exchanger's design parameters and the conditions of the hot storage. The second function is to compute the pressure drop and friction losses across the exchanger shells resulting from the molten salt flow. This is done using equations proposed by Gaddis and Gnielinsky~\cite{gaddis1997pressure}.
Another way to model the heat exchange to produce superheated steam from molten salt is to suppose a perfect heat exchanger. Both approaches can be selected separately on different \solar{} instances. In the idealized version, the heat exchanger has a 100\% thermal efficiency without friction loss and the parameters describing its geometry are not required. 

\subsection{Molten salt cycle}


Four subsystems are part of the molten salt cycle: the receiver, the two storage tanks, and the heat exchanger. The state and requirements of these subsystems depends on the state of other systems. For instance, the molten salt rate required by the heat exchanger to generate the necessary amount of steam at any moment depends not only on its own design parameters but also on the state of the hot storage: a higher inlet temperature will result in a smaller molten salt flow. Though because the rate of temperature drop of a storage unit is dependent on the stored mass, the two components are mutually dependent.


Similarly, the amount of molten salt that can be taken from the cold storage and heated to the design point conditions through the receiver is dependent on the temperature of the cold storage. Needless to say that the outlet conditions of both the receiver and the exchanger will too impact the storage units temperature and levels.

\subsection{Power block}

The only  variable related to the power block is the choice of the turbine  model. Technical data for a variety of  commercial steam turbines used  in CSP applications is retrieved from
Siemens~\cite{steamturbinessiemens}.  The dataset includes steam inlet pressures and temperatures required to operate each turbine, as well as their respective maximum and minimum power output. This part of the model is used to determine the amount of  thermal power that ought to be extracted from the molten salt in order to meet the power demand.

Turbine efficiency is computed as a function of the inlet steam conditions and capacity usage ratio using a simple empirical model provided by Bahadori and Vuthaluru~\cite{Bahadori2010}.
Steady-state conditions at all times are assumed, and turbine transients are neglected. A mechanical-to-electrical efficiency of 95\% is assumed.

All components of the power cycle other than the turbine are idealized. The power required  to pump the water condensate towards the heat exchanger is treated as part of the parasitic power consumption. No transient regime is considered and the whole cycle is assumed to shift instantly to match the demand profile. For most turbines, there exists a minimum value of power for which  they can be operated. In the event that the demand is inferior to the minimum requirement for a turbine, the model will operate it at its minimum, least efficient regime, if possible, in order to reflect the fact that stopping the plant entirely in the middle of production is usually a bad operational option.

\subsection{Auxiliary constraints models}

Auxiliary models are used to provide constraints to the optimization problems.
Without them, the solution to many of the problems would be trivial  or impractical. 
For example, if the  thickness of the insulation on the tanks is not counter-balanced by any  other factor,  it will always  be set at  its highest value. 
Similar reasoning applies for the heliostats: 
	 the best way to maximize the field's surface efficiency is to fill it with a maximum number of very small heliostats that cover the entire field.

In order to provide a sufficient amount of constraints, four auxiliary models  are used to estimate the equipment cost, parasitic loads, and the pressure  both in the tubes of the receiver and heat exchanger.

\subsubsection{Initial capital cost model}

Although no complete life cycle cost analysis is integrated in the simulation, a simple initial investment cost model is provided in order to serve as a limiting factor for many of the design parameters. The data used to build the capital cost model was taken mostly from the National Renewable Energy Laboratory report on the SAM (System Advisor Model)~\cite{turchi2013molten} software and the Sandia Roadmap report~\cite{kolb2011power}. The SAM model provides a means to determine the value of each component, but does so through an empirical model that mostly uses relations to the size or desired capacity of the plant. Its objective is to predict the potential cost of a project based on its sheer scale, rather than on specific technical characteristics. This turns out to be of little use as a limiting factor in optimization problems.

In order to link the cost of the components to their respective design parameters, relations have been established based on the price of materials. 
The total cost of the power plant is given below and is the sum of the costs of its subsystems, each  one described as a function of some of the design parameters previously listed in their respective sections.
These parameters also appear as the optimization variables of Tables~\ref{tab-vars1} and~\ref{tab-vars2}.

\begin{equation}
\begin{split}
\mbox{Total Cost} = \ 
& 
\underbrace{
    N_{\hf}
    C_{\hf} ( L_{\hf},  W_{\hf})
}_\text{Heliostats field} +
\underbrace{
    C_{\twr}(H_{\twr})
}_\text{ Tower} + 
\underbrace{
    C_{\rcv}(H_{\rcv},W_{\rcv},t_{\rcv})}_\text{Receiver} + \\
&
\underbrace{
    C_{\hot}(H_{\hot}, d_{\hot}, t_{\hot})
}_\text{ Hot storage} +
\underbrace{
    C_{\cold}(H_{\cold}=H_{\hot}\times 1.2, d_{\cold}=d_{\hot}, t_{\cold})
}_\text{ Cold storage} + \\
&
\underbrace{
    C_{\sg}(S_t, L_{\sg}, d_{\sg}, D_{\sg}, H_{\sg,\baf}, N_{\sg,\baf}, N_{\sg,\tub}, N_{\sg,\sh,\pass}, N_{\sg,\tub,\pass})   
}_\text{ Steam generator} + \\
&
  \underbrace{C_{\trb} (ST)}_\text{ Turbine}.
\end{split}
\label{eq:totalCost}
\end{equation}

\subsubsection{Parasitic loads model}

The parasitic  load is the power required to operate the plant. The SAM considers a detailed set of parasitic  load sources, including the  heliostats electronic and  mechanical control systems, piping anti-freeze protections, the  molten salt pumps in the storage units,  and the  control systems of the receiver, heat exchanger  and steam condenser. In the current work, components that  are not directly necessary to describe the  operation of a CSP  plant  (e.g., pump electronic controls) have been idealized. Since the pumps and pipes linking the main components of the system are not explicitly simulated or subject to optimization, their contribution to the losses  is not considered. The evaluation of the parasitic load  $P$ is thus comprised of five terms:

\begin{equation*}
P = 
\underbrace{P_{\rcv}}
_\text{Receiver} + 
\underbrace{P_{\sg,\sh}}
_{\substack{\text{Steam}\\ \text{generator} \\ \text{ (shell)}}} + 
\underbrace{P_{\sg,\tub}}
_{\substack{\text{Steam}\\ \text{generator} \\ \text{(tubes)}}} + 
\underbrace{P_{\hf}}
_{\substack{\text{Heliostats} \\ \text{field} }} + 
\underbrace{P_{\hot} + P_{\cold}}
_{\substack{ \text{Storage} \\ \text{anti-freezing} }}.
\label{parasitic_gen}
\end{equation*}

The values of  $P_{\rcv}$ and  $P_{\sg,\tub}$ are computed as the power required to pump a fluid through smooth tubes assuming a constant pump efficiency of 90\%.
For  $P_{\sg,\sh}$,
a shell-side pressure drop model from~\cite{gaddis1997pressure} is used for shell-and-tubes heat exchangers.
It is assumed that the pumps are capable of providing the required flow rate and pressure. The pressure differential is a function of both the flow rate and the geometry of the components. The term  $P_{\hf}$ considers a constant and identical power  consumption of  55~$\W$  per heliostat  when operated (i.e., daytime)~\cite{amadei2013simulation}.  When the storage temperature drops below a certain threshold, the tanks are heated with  $P_{\hot} + P_{\cold}$ to keep the temperature of both storage units above the molten salt's melting point.

\subsubsection{Maximum allowable pressure in tubes}

 The thickness and diameter of the tubes in both the receiver and steam generator is limited by the stress imposed on the steel  by the internal fluid pressure. The stress  is checked in each tube cross section  so that it never exceeds the steel's yield strength. For the sake of simplicity and time, creeping effects in metals are not considered.  However, given that the system operates at  high temperatures, this effect  should be accounted for in more accurate models.

An exhaustive study of the heat exchanger system would also consider the stress imposed on the tubes by the molten salt flowing across them, thereby providing an additional limiting factor on the tubes length and spacing and the baffles spacing. In the present work, though, only the outward radial pressure exerted by the fluid flowing inside the tubes is considered.

\section{A collection of optimization problem instances}
\label{sec-instances}
The previous sections described the building blocks necessary to construct 
    a collection of 10 optimization problem instances,
    provided by the \solar{} code.
This collection presents different variations of BBO benchmarks which can be used to test a wide range of algorithms. Each instance is uniquely characterised by its set of variables (Section~\ref{sec-vars}), its objectives (Section~\ref{sec-obj}), its constraints (Section~\ref{sec-cstr}) and its chosen fixed parameters, each of which is described in this section.

\subsection{A general class of blackbox optimization  problems}

The instances of the collection belong to the following general minimization problem:

\begin{align*}
\min\limits_{\xx \in \R^n} \ F(\xx)=\bigl(f_1(\xx), f_2(\xx), \ldots,f_p(\xx)\bigr)
\end{align*}
\begin{subnumcases}{\textrm{subject to } ~~  \label{eq-pb}}
    c_i(\xx) \leq 0     &  $i\in\{1,2,\ldots,m\}$ \label{eq-pb1} \\
    x_j \in \Z        &  $j \in {\cal J}$ \label{eq-pb2} \\
    \ell \leq \xx \leq u                    \label{eq-pb3}
\end{subnumcases}
%
where
$F:\R^n \rightarrow \RU^p$  is the objective function(s) to minimize,
with $\RU := \R \cup \{ \pm \infty\}$,
subject to Constraints~\eqref{eq-pb1},~\eqref{eq-pb2}, and~\eqref{eq-pb3},
and ${\cal J} \subseteq \{1,2,\ldots,n\}$ is the set of discrete variables.
Constraint~\eqref{eq-pb1} 
    is compactly written as $C(\xx) \leq 0$ where
    $C(\xx)=\bigl( c_1(\xx), c_2(\xx), \ldots, c_m(\xx) \bigr)$
    and with
    $C:\R^n \rightarrow \RU^m$.
The vector $\xx \in \R^n$ is the vector of optimization variables  with bounds
defined by~\eqref{eq-pb3}.

In practice, each component of $\xx$ may be real, integer, binary or categorical.
The indices of the integer and binary variables belong to the set ${\cal J}$.
Categorical variables are also represented by integers. 
    

While the underlying  model is known and described in the previous section, 
	the actual implementation is treated as a blackbox. 
For any value of the variable $\xx$ and for any instance with
 	$p$ objective functions and $m$ constraints, 
	the computer simulation returns a formatted output in the form of a single vector $y$ 
	containing the values of $F(\xx)$ and $C(\xx)$ 
	such that 
\begin{equation}
y(\xx) \ = \ \bigl( F(\xx),C(\xx) \bigr)
    \ =\ \bigl( f_1(\xx), f_2(\xx), \ldots, f_p(\xx), c_1(\xx), c_2(\xx), \ldots, c_m(\xx) \bigr)
\label{eq-outputs}
\end{equation}

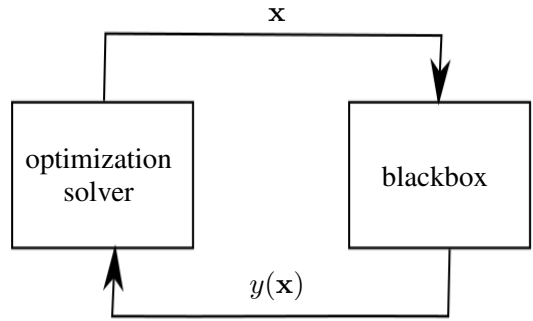
\begin{figure}[ht!]
\setlength{\unitlength}{10mm}
\centering
\includegraphics[width=7cm, height=4cm]{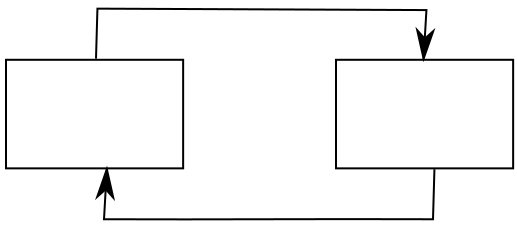}
\begin{picture}(0,0)(7,0)
\put(5.6,1.95){\makebox(0,0){blackbox}}
\put(1.13,2.15){\makebox(0,0){optimization}}
\put(1.13,1.75){\makebox(0,0){solver}}
\put(3.5,4.1){\makebox(0,0){$\xx$}}
\put(3.5,0.5){\makebox(0,0){$y(\xx)$}}
\end{picture}
\caption{Blackbox optimization (BBO) process.}
\label{fig-OptProcess}
\end{figure}

The optimization process illustrated in Figure~\ref{fig-OptProcess}
	consists in successively calling the blackbox with different values of $\xx$, 
	and using the corresponding output values $y(\xx)$ 
    to dictate the next input values provided to the blackbox.

\subsection{Ten blackbox instances}

Ten instances, named~\solar{1.1} to~\solar{10.1}, are implemented in a single {\sf C++} program
that can be executed as a console application.
The ``{\sf .1}'' notation refers to the release {\sf 1.0} of the code, which is now frozen in time and guaranteed to never change.
If bugs that may affect the outputs of the simulator are discovered, they will be fixed in another release ({\sf 2.0} or more).
This is the reason why users of \solar{} are strongly encouraged to use the notation \solar{X.1}.

In order to compute $y(\xx)$, 
    the following elements can be supplied to the \solar{} code:
\begin{itemize}
    \item the instance number (also called ``problem~ID'' in the code) is an integer ranging from 1 to 10 that indicates which one of the 10 instances should be used;
    \item the input vector $\xx$;
    \item the ``{\tt seed}'' value, for the random seed of stochastic instances (see Section~\ref{subsec-stoch});
    \item the ``{\tt rep}'' value, for the number of replications for the stochastic instances (see Section~\ref{subsec-stoch});
    \item the ``{\tt fid}'' value, to select the fidelity of the output (see Section~\ref{subsec-multi-fid}). 
\end{itemize}

Only the instance number and the vector $\xx$ are mandatory. The default values for the remaining elements correspond to deterministic instances with the best fidelity.


Each one of the 10 instances is characterized by a specific format for 
	$\xx$ (number of variables and their respective type and bounds) and 
	$y$ (the objective function(s) and constraints),
    as summarized in Table~\ref{tab-characteristics}. 

\begin{table}[ht!]
\renewcommand{\tabcolsep}{4.5pt}
\begin{footnotesize}
\begin{center}
\begin{tabular}{|l|cc|c|c|cc|c|c|c|}
\hline
Instance & \multicolumn{3}{c|}{\# of variables} & \# of obj.  &
           \multicolumn{3}{c|}{\# of constraints}   & \# of stoch. outputs &
           Multi    \\
         & cont. & discr. (cat.) & $n$  &  $p$
         & simu. & a~priori (lin.) & $m$ & (obj. or constr.)     & fidelity \\
\hline
\hline
\solar{1.1}  &  8 & 1\quad (0) &  9 & 1 &  2 & 3\quad (2) &  5 & 1 &  no \\
\solar{2.1}  & 12 & 2\quad (0) & 14 & 1 &  7 & 5\quad (3) & 12 & 4 & yes \\
\solar{3.1}  & 17 & 3\quad (1) & 20 & 1 &  8 & 5\quad (3) & 13 & 5 & yes \\
\solar{4.1}  & 22 & 7\quad (1) & 29 & 1 &  9 & 7\quad (5) & 16 & 6 & yes \\
\solar{5.1}  & 14 & 6\quad (1) & 20 & 1 &  8 & 4\quad (3) & 12 & 0 &  no \\
\solar{6.1}  &  5 & 0\quad (0) &  5 & 1 &  6 & 0\quad (0) &  6 & 0 &  no \\
\solar{7.1}  &  6 & 1\quad (0) &  7 & 1 &  4 & 2\quad (1) &  6 & 3 & yes \\
\solar{8.1}  & 11 & 2\quad (0) & 13 & 2 &  4 & 5\quad (3) &  9 & 4 & yes \\
\solar{9.1}  & 22 & 7\quad (1) & 29 & 2 & 10 & 7\quad (5) & 17 & 6 & yes \\
\solar{10.1} &  5 & 0\quad (0) &  5 & 1 &  0 & 0\quad (0) &  0 & 0 &  no \\
\hline
\end{tabular}
\end{center}
\end{footnotesize}
\caption{Characteristics of the ten blackbox problem instances, including the number of continuous and discrete variables (including categorical), the number of objectives, the number of simulation and a~priori constraints (see~\cite{LedWild2015} and Section~\ref{subsec-cstr-eval}), and the number of stochastic outputs available through the {\tt -seed} option (see Section~\ref{subsec-stoch}).
The last column indicates if the multifidelity feature is available through the {\tt -fid} option (see Section~\ref{subsec-multi-fid}).
\solar{10.1} is the unconstrained version based on \solar{6.1}.
\solar{2.1} is the only instance with a known analytical expression for the objective function.
}
\label{tab-characteristics}
\end{table} 
 

\subsubsection{Objectives and general description of each instance}
\label{sec-obj}

This section gives a high-level description of each problem instance,
and in particular, it describes the objective functions to optimize.
When $p=1$, this is single-objective optimization problem, and when $p=2$, this is a biobjective optimization problem.
The optimization is subject to various types of constraints that are described in Section~\ref{sec-cstr}.

\paragraph{\solar{1.1} - Maximize heliostats field energy output}
This instance runs only the heliostats field model. It uses 9 variables of which one is discrete and the others are continuous. The objective is to maximize the  energy collected by the receiver in 24 hours, while respecting a \$50M budget and a maximum field area. The objective is subject to 5 relaxable and quantifiable constraints.
 
\paragraph{\solar{2.1} - Minimize the heliostats field surface (analytical objective)}
This instance runs the whole power plant model and uses the idealized model for the heat exchanger. It uses 14 variables of which 2 are discrete and the others are continuous. The objective is $F(\xx)=f_1(\xx)={x_3^2(x_9^2-x_8^2) x_7}\frac{\pi}{180}$
which corresponds to the heliostats field surface to minimize while satisfying the power demand peaking at 20~M$\W$ and respecting a \$300M budget. The objective is subject to 14 relaxable and quantifiable constraints.

\paragraph{\solar{3.1} - Minimize total investment cost}
This instance simulates the whole power plant and uses the idealized model for the steam generator. It uses 20 variables of which two are discrete, one is a categorical variable and 17 are continuous. The objective is to minimize the total investment cost while satisfying the demand and respecting a maximum field size. The simulation is done over 24 hours and the power plant is required to provide constant 10~M$\W$ during peak hours, between noon and 6~pm The objective is subject to 12 relaxable and quantifiable constraints.

\paragraph{\solar{4.1} - Minimize total investment cost}
It is similar to \solar{3.1}, but with an increased level of complexity in the NTU-effectiveness steam generator model presented in Section~\ref{sec-heat-exchanger}.
It uses 29 variables of which 6 are discrete, one is a categorical variable and 22 are continuous. The objective is to minimize the total investment cost while satisfying the demand and respecting a maximum field size. The simulation is done over 72 hours and the power plant is required to provide power at all time for a summer day demand profile peaking at 25~M$\W$. The objective is subject to 16  relaxable and quantifiable constraints.

\paragraph{\solar{5.1} - Maximize the satisfaction of the demand}
This instance runs the  molten salt loop and  uses estimates of the performance of the heliostats field from a pre-computed database in order to reduce the computation time. It uses 20 variables of which 5 are discrete, one is a categorical variable and 14 are continuous. The power plant performance is simulated over a period of 30 days with an inconsistent field performance analogous to slightly unreliable weather conditions. The objective is to maximize the time for which the power plant is able to operate at nominal capacity. A surrogate model is available for this problem which consists in running the simulation on only a fraction of the 30 days and extrapolating the resulting performance over the 30 days. The objective is subject to 12  relaxable and quantifiable constraints.

\paragraph{\solar{6.1} - Minimize the cost of storage}
This instance runs a predetermined power plant using the  molten salt cycle and power block models. It uses 5 continuous variables. The objective is to minimize the cost of the thermal storage units so that the power plant is able to sustain a 100~M$\W$ electrical power output  during 24 hours. Since the heliostats field is not being optimized, its hourly power output is read from a prerecorded file instead of being simulated,  thus reducing the computation time.
The objective is subject to 2 relaxable and quantifiable constraints and 4 unrelaxable and quantifiable constraints.

\paragraph{\solar{7.1} - Maximize receiver efficiency}
This instance simulates the heliostats field and the receiver unit over a 24  hour period. It uses 7 variables, of which one is discrete and the others are continuous. The objective is to maximize the receiver efficiency. A surrogate version of the model can be used so that a much lower density of sunrays is used to evaluate the field's performance. The objective is subject to 6 unrelaxable and quantifiable constraints.

\paragraph{\solar{8.1} - Maximize heliostats field performance and minimize cost (biobjective)}
This instance runs the heliostats field and receiver models. It uses 13 variables, of which two are discrete and 11 are continuous. This is a biobjective problem of which the two objectives are to maximize the amount of energy transferred to the molten salt over a 24 hours period all while minimizing the total cost of the field, tower and receiver. The optimization is conducted over the design parameters of both the heliostats field and the receiver. The objectives are subject to 9 relaxable and quantifiable constraints.

\paragraph{\solar{9.1} - Maximize power and minimize losses (biobjective)}
This instance simulates the entire power plant over a single day. It uses 29 variables, of which 6 are discrete, one is a categorical variable and 22 are continuous. This is a biobjective problem of which the two objectives are to maximize the generated electrical power and minimize the parasitic losses while respecting a \$1.2B budget. The objectives are subject to 5   relaxable and quantifiable constraints and 12 unrelaxable and quantifiable constraints.

\paragraph{\solar{10.1} - Minimize the cost of storage (unconstrained)}
This instance is the unconstrained version of \solar{6.1} where the constraints are penalized in the objective with
$$F(\xx)=f_1(\xx)= \frac{h_1(\xx)}{10^6} + \frac{\left( g_1(\xx)^2 + (2 \times 10^{-6} g_2(\xx))^2 + g_3(\xx)^2 + g_4(\xx)^2 + g_5(\xx)^2 + g_6(\xx)^2 \right)}{2}$$
with $h_1$ the objective function of \solar{6.1} and 
$g_i$ the violation of Constraint~$i$  of \solar{6.1} with $i\in\{1,2,\ldots,6\}$.
These weights have been fixed empirically.

\newpage

\begin{table}[htb!]
\renewcommand{\tabcolsep}{1.25pt}
\centering
{\footnotesize
\begin{tabular}{|c|clc|crr|c|c|}
\hline
\multirow{2}{*}{Category} & \multirow{2}{*}{Symbol} & \multirow{2}{*}{Quantity} 	&
\multirow{2}{*}{Unit}&\multirow{2}{*}{Type}	& {Lower} 	& {Upper} 	
 	& \multirow{2}{*}{Instances} & Optim.	\\
&  						&  				&   		&  		&	 {bound}	& { bound}	& 	&  variable				\\
\hline
\hline
\multirow{11}{*}{Heliostats field}
& $L_{\hf}$ &   Heliostats length		 			& $\m$ & cont.  			& 1 			& 40 					&1~2~3~4~8~9& $x_1$	\\
\cline{2-9}
&  $W_{\hf}$ & Heliostats width  		 			& $\m$ & cont. 	 			& 1 			& 40 					&1~2~3~4~8~9	& $x_2$	\\
\cline{2-9}
& $H_{\twr}$   & Tower height             		    		& $\m$ & cont. 			& 20 			& 250 			&1~2~3~4~8~9 & $x_3$		\\
\cline{2-9}
& \multirow{2}{*}{$H_{\rcv}$} &\multirow{2}{*}{Receiver aperture height}	 			& \multirow{2}{*}{$\m$} 	& \multirow{2}{*}{cont.} 			& \multirow{2}{*}{1} 			& \multirow{2}{*}{30} 			 		&1~2~3~4~8~9	& $x_4$	\\
&&&&&&      &7	& $x_1$	\\
\cline{2-9}
	& \multirow{2}{*}{$W_{\rcv}$} &\multirow{2}{*}{Receiver aperture width}   			& \multirow{2}{*}{$\m$} 	& \multirow{2}{*}{cont.} 			& \multirow{2}{*}{1} 			& \multirow{2}{*}{30} 					&1~2~3~4~8~9	& $x_5$	\\
&&&&&&        &7	& $x_2$	\\
\cline{2-9}
& $N_{\hf}$& Number of heliostats to fit    	  		& $-$ 	& 	discr.	& 1  			& $+\infty$	& 1~2~3~4~8~9& $x_6$		\\
\cline{2-9}
& $\theta_{\hf}$ & {Field angular width}  & $\deg$ & cont.& 1 & 89 & 1~2~3~4~8~9 & $x_7$ \\
\cline{2-9}
& $R_{\hf}^{\min}$ 	& Min. distance from tower &  ratio	& cont. & 0 & 20 & 1~2~3~4~8~9 & $x_8$\\
\cline{2-9}
& $R_{\hf}^{\max}$ 	& Max. distance from tower &  ratio	& cont. & 1 & 20 & 1~2~3~4~8~9 & $x_9$\\
\hline
\hline
& \multirow{3}{*}{$T_{\rcv}^{out}$} & \multirow{3}{*}{Receiver outlet temp.} & \multirow{3}{*}{$\K$} & \multirow{3}{*}{cont.} & \multirow{3}{*}{793} & \multirow{3}{*}{995} & 2~3~4~9 & $x_{10}$ \\
&  &  & &  & &                                             & 5~6~10         & $x_1$\\
&  &  & &  & &                                             & 7             & $x_3$ \\
\cline{2-9}
 & \multirow{3}{*}{$H_{\hot}$} & \multirow{3}{*}{Hot storage height} & \multirow{3}{*}{$\m$} & \multirow{3}{*}{cont.}
            & 1 & 50 &  3~4~9 & $x_{11}$ \\	
 &&&& &  1 & 30 & 5 & $x_2$ \\
 &&&& &  2 & 50  &6~10  & $x_2$ \\
\cline{2-9}
& \multirow{3}{*}{$d_{\hot}$}  & \multirow{3}{*}{Hot storage diameter} & \multirow{3}{*}{$\m$} & \multirow{3}{*}{cont.}
            & 1 & 30 & 3~4~9 & $x_{12}$ \\
 &&&& & 1 & 30 & 5         & $x_3$ \\
 &&&& & 2 & 30 & 6~10         & $x_3$ \\
\cline{2-9} 
 & \multirow{3}{*}{$t_{\hot}$} & \multirow{3}{*}{Hot storage insulation thickness} & \multirow{3}{*}{$\m$} & \multirow{3}{*}{cont.}
           & 0.01 & 5  & 3~4~9 & $x_{13}$ \\ 
 &&&& & 0.01 & 2  & 5 & $x_4$\\
 &&&& & 0.01 & 5  & 6~10 & $x_4$ \\
 \cline{2-9}
Heat transfer loop
& \multirow{3}{*}{$t_{\cold}$} & \multirow{3}{*}{Cold storage insulation thickness} & \multirow{3}{*}{$\m$} & \multirow{3}{*}{cont.}
           & 0.01 & 5  & 3~4~9 & $x_{14}$ \\ 
 &&&& & 0.01 & 2  & 5 & $x_5$\\
 &&&& & 0.01 & 5  & 6~10 & $x_5$ \\
\cline{2-9}
& \multirow{2}{*}{$T_{\cold}^{\min}$} & \multirow{2}{*}{Min. cold storage temp.} & \multirow{2}{*}{$\K$} & \multirow{2}{*}{cont.} & \multirow{2}{*}{495} & \multirow{2}{*}{650}
                    & 3~4~9 & $x_{15}$ \\
 &&&& &  &  & 5 & $x_{6}$ \\
 \cline{2-9}
& \multirow{5}{*}{$N_{\rcv,\tub}$} &  \multirow{5}{*}{Receiver number of tubes} & \multirow{5}{*}{$-$} & \multirow{5}{*}{discr.} & \multirow{5}{*}{1}            & 9,424  & 2& $x_{11}$ \\
&  &   & &  &  & 9,424  & 3 & $x_{16}$ \\
&  &   & &  &  & 7,853   & 4~9& $x_{16}$ \\
&  &   & &  &  & 1,884  & 5& $x_{7}$ \\
&  &   & &  &  & 8,567 & 7& $x_{4}$ \\
&  &   & &  &  & 7,853 & 8& $x_{10}$ \\
\cline{2-9}
& \multirow{5}{*}{$t_{\rcv}$} & \multirow{5}{*}{Receiver insulation thickness} & \multirow{5}{*}{$\m$} & \multirow{5}{*}{cont.}
          & 0.01 &5  & 2 	 &   $x_{12}$ \\
&&&& & 0.01 &5  & 3~4~9	 & $x_{17}$ \\
&&&& & 0.10 & 2  & 5	&   $x_8$ \\
&&&& & 0.01 & 5 & 7	&   $x_5$ \\
(continued in Table~\ref{tab-vars2}) &&&& & 0.01 & 5  & 8	&    $x_{11}$ \\
\hline
\end{tabular}
}
\caption{The 29 possible optimization variables (1/2).
The ``$\times H_{\twr}$'' unit means that the distances from tower ($x_8$ and~$x_9$) are expressed as multiples of 
$H_{\twr}$.}
\label{tab-vars1}
\end{table}

\subsubsection{Optimization variables}
\label{sec-vars}

The optimization variables of Problem~\eqref{eq-pb} are $\xx \in \R^n$.
While all the variables are by definition continuous, Equation~\eqref{eq-pb2} imposes
that $x_j \in \Z$ for $j \in {\cal J}$.
In addition, the type of turbine ($ST$) is a categorical variable with 8 possible values.
Tables~\ref{tab-vars1} and~\ref{tab-vars2} list all possible 29 variables, their characteristics and bounds ($\ell \leq \xx \leq u$~\eqref{eq-pb3}),
depending in the instance.
Only Instances~\solar{4.1} and~\solar{9.1} include all of them as optimization variables.
Some optimization variables have no upper bounds. As some algorithms require all bounds to be defined, this task is left to the user that can create a reasonable upper bound from the suggested lower bounds and initial values.
The simulation parameters and the values of the instance-specific constraints are chosen to reflect realistic engineering problems. For instance, budgets, maximum field surface, and power output requirements are derived from known data on existing similar power plants.
The exact numerical values to which these numerous parameters are set in each problem instance
 can be found in~\cite{MScMLG}.


 

\begin{table}[ht!]
\renewcommand{\tabcolsep}{1pt}
\centering
{\footnotesize
\begin{tabular}{|c|clc|crr|c|c|}
\hline
\multirow{2}{*}{Category} & \multirow{2}{*}{Symbol} & \multirow{2}{*}{Quantity} 	&
\multirow{2}{*}{Unit}&\multirow{2}{*}{Type}	& {Lower} 	& {Upper} 	
 	& \multirow{2}{*}{Instances} & Optim.	\\
&  						&  				&   		&  		&	 {bound}	& { bound}	& 	&  variable				\\
\hline
\hline
& \multirow{5}{*}{$d_{\rcv}$} & \multirow{5}{*}{Receiver tubes inner diameter} & \multirow{5}{*}{$\m$} & \multirow{5}{*}{cont.}
          &  \multirow{5}{*}{0.005}&\multirow{5}{*}{0.1} & 2 & $x_{13}$ \\
&&&& & & & 3~4~9& $x_{18}$ \\
&&&& && & 5& $x_{9}$ \\
&&&& && & 7& $x_{6}$ \\
Heat transfer loop &&&& && & 8& $x_{12}$ \\
\cline{2-9}
(continued from Table~\ref{tab-vars1})
& \multirow{6}{*}{$D_{\rcv}$} & \multirow{6}{*}{Receiver tubes outer diameter} & \multirow{6}{*}{$\m$} & \multirow{6}{*}{cont.}
	   & 0.0050&0.1 & 2 & $x_{14}$ \\
&&&& & 0.0050&0.1 & 3 & $x_{19}$ \\
&&&& &0.0060&0.1 & 4~9 & $x_{19}$ \\
&&&& &0.0050&0.1 & 5 & $x_{10}$ \\
&&&& &0.0055&0.1 & 7 & $x_{7}$ \\
&&&& &0.0060&0.1 & 8 & $x_{13}$ \\
\hline
\hline
\multirow{18}{*}{Steam generator}
& \multirow{2}{*}{$S_t$} & \multirow{2}{*}{Tubes spacing} & \multirow{2}{*}{$\m$} & \multirow{2}{*}{cont.}
                & 0.007 & 0.2 & 4~9 & $x_{20}$ \\
&  &  &  & & 0.006 & 0.2 & 5 & $x_{11}$\\
\cline{2-9}
& \multirow{2}{*}{$L_{\sg}$} & \multirow{2}{*}{Tubes length} & \multirow{2}{*}{$\m$} & \multirow{2}{*}{cont.} & \multirow{2}{*}{0.5} & \multirow{2}{*}{10}
                        & 4~9 & $x_{21}$\\
&  & &  &  &  & & 5 & $x_{12}$\\
\cline{2-9}
& \multirow{2}{*}{$d_{\sg}$} & \multirow{2}{*}{Tubes inner diameter} & \multirow{2}{*}{$\m$} & \multirow{2}{*}{cont.} & \multirow{2}{*}{0.005} &\multirow{2}{*}{0.1} &4~9 & $x_{22}$\\
&  &  &  &  &  & &5 & $x_{13}$\\
\cline{2-9}
&\multirow{2}{*}{$D_{\sg}$} & \multirow{2}{*}{Tubes outer diameter} & \multirow{2}{*}{$\m$} & \multirow{2}{*}{cont.} &\multirow{2}{*}{0.006} & \multirow{2}{*}{0.1} &4~9 & $x_{23}$\\
&  &  &  &  &  & &5 & $x_{14}$\\
\cline{2-9}
& \multirow{2}{*}{$H_{\sg,\baf}$} & \multirow{2}{*}{Baffles cut} & \multirow{2}{*}{ratio} & \multirow{2}{*}{cont.} & \multirow{2}{*}{0.15} & \multirow{2}{*}{0.4} &4~9 & $x_{24}$\\
&  &  &  &  &  & &5 & $x_{15}$\\
\cline{2-9}
(exchanger)
& \multirow{2}{*}{$N_{\sg,\baf}$} &\multirow{2}{*}{Number of baffles} & \multirow{2}{*}{$-$} & \multirow{2}{*}{discr.} & \multirow{2}{*}{2} & \multirow{2}{*}{$+\infty$} & 4~9 & $x_{25}$\\
&  &  &  &  &  & &5 & $x_{16}$\\
\cline{2-9}
& \multirow{2}{*}{$N_{\sg,\tub}$} & \multirow{2}{*}{Number of tubes} &\multirow{2}{*}{$-$} & \multirow{2}{*}{discr.} & \multirow{2}{*}{1} & \multirow{2}{*}{$+\infty$} &4~9 & $x_{26}$\\
&  &  &  &  &  & &5 & $x_{17}$\\
\cline{2-9}
& \multirow{2}{*}{$N_{\sg,\sh,\pass}$} & \multirow{2}{*}{Number of shell passes} &\multirow{2}{*}{$-$}& \multirow{2}{*}{discr.} & \multirow{2}{*}{1} & \multirow{2}{*}{10} &4~9 & $x_{27}$\\
&  &  &  &  &  & &5 & $x_{18}$\\
\cline{2-9}
& \multirow{2}{*}{$N_{\sg,\tub,\pass}$} & \multirow{2}{*}{Number of tube passes} &\multirow{2}{*}{$-$}& \multirow{2}{*}{discr.} &  \multirow{2}{*}{1} & \multirow{2}{*}{9} &4~9 & $x_{28}$\\
&  &  &  &  &  & &5 & $x_{19}$\\
\hline
\hline
\multirow{2}{*}{Powerblock} 
& \multirow{2}{*}{$ST$} & \multirow{2}{*}{Type of turbine} & \multirow{2}{*}{$-$} & discr. & \multirow{2}{*}{1} & \multirow{2}{*}{8} & 3~5 & $x_{20}$ \\
&  &  & & (cat.) &  &  & 4~9 & $x_{29}$ \\
\hline
\end{tabular}
}
\caption{The 29 possible optimization variables (2/2).}
\label{tab-vars2}
\end{table}

\subsubsection{Constraints}
\label{sec-cstr}

The taxonomy of constraints of~\cite{LedWild2015} allows to better describe the different types of constraints in Problem~\eqref{eq-pb}.
Constraints~\eqref{eq-pb1} and~\eqref{eq-pb3} (bounds) are {\em quantifiable}: the inequality format ensures that $c_i(\xx)$ ($i \in \{1,2,\ldots,m\}$) or $\xx$ measure the distance to feasibility or infeasibility. On the contrary, Constraints~\eqref{eq-pb2} (discrete nature of some variables) are {\em nonquantifiable}.
Constraints~\eqref{eq-pb2} and~\eqref{eq-pb3} (discrete variables and bounds) are also {\em a~priori} constraints: they can be evaluated without executing the \solar{} simulator.
Depending on the instances, some of the constraints of~\eqref{eq-pb1} are also a~priori (including linear constraints), while the others are of the type {\em simulation} (the execution of the \solar{} simulator is necessary to evaluate them).
A~priori constraints are also considered to be {\em unrelaxable}:
They represent known relations between some variables and need to always be respected for the simulator to be executed: for example, the inner diameter of a tube must be smaller than its outer diameter. 
The simulation constraints, on the contrary, are {\em relaxable}, meaning that their potential violation does not prevent to execute the simulation.
Hence, they can be violated during an optimization process, which can be very useful,
as long as the final proposed solution is feasible.
%
%
Finally, there are {\em hidden} constraints consisting mostly of code glitches or instability; 
	some solutions may cause the code to crash or  otherwise malfunction. 
These constraints are not given explicitly in the problem definition as they are unknown.
%
Section~\ref{subsec-cstr-eval} provides more details on how \solar{} manages a~priori and hidden constraints.


The \solar{} code, for a given $\xx$, returns $y(\xx)$ as defined in~\eqref{eq-outputs},
which, for the constraints, corresponds to the left hand-sides $c_i(\xx)$ of Constraints~\eqref{eq-pb1},
for $i \in \{1, 2,\ldots, m\}$.
Constraints~\eqref{eq-pb2} and~\eqref{eq-pb3} must be considered by the optimization method,
and if they are not satisfied, the \solar{} simulator will not execute.
The bounds and the discrete variables are visible in Tables~\ref{tab-vars1} and~\ref{tab-vars2}.
Hence the following of this section focuses on how Constraints~\eqref{eq-pb1} are defined in
\solar{1.1} to~\solar{9.1} (\solar{10.1} is unconstrained).
The number and type of these constraints are given in Table~\ref{tab-characteristics}.

The a~priori constraints are given below, using the $\xx$ vector. However, since $\xx$ is not
defined the same for each instance, it is necessary to consult Tables~\ref{tab-vars1} and~\ref{tab-vars2}.


\subparagraph{Linear a~priori constraints:}

Constraints details are shown in Table~\ref{tab-linear-apriori-constraints}.

\begin{table}[H]
\renewcommand{\tabcolsep}{4.5pt}
\begin{footnotesize}
\begin{center}
\begin{tabular}{|l|ccccccccc|}
\hline
& \solar{1.1} & \solar{2.1} & \solar{3.1} & \solar{4.1} & \solar{5.1} & \solar{6.1} & \solar{7.1} & \solar{8.1} & \solar{9.1} \\
\hline
\hline
Constraint & \multicolumn{9}{c|}{Tower is at least twice as high as heliostats} \\
Indexation & $c_3$ & $c_4$ & $c_3$ & $c_3$ & $-$ & $-$ & $-$ & $c_2$ & $c_4$ \\
[.2cm]
Equation & \multicolumn{9}{c|}{$2x_1 \leq x_3$} \\
[.2cm]
\hline
\hline
Constraint & \multicolumn{9}{c|}{Minimum distance from tower is lower than the maximum distance from tower} \\
Indexation & $c_4$ & $c_5$ & $c_4$ & $c_4$ & $-$ & $-$ & $-$ & $c_3$ & $c_5$ \\
[.2cm]
Equation & \multicolumn{9}{c|}{$x_8 \leq x_9$} \\
[.2cm]
\hline
\hline
Constraint & \multicolumn{9}{c|}{Receiver tubes inside diameter is lower than the outside diameter} \\
Indexation & $-$ & $c_{11}$ & $c_{10}$ & $c_{10}$ & $c_6$ & $-$ & $c_3$ & $c_6$ & $c_{11}$ \\
[.2cm]
Equation & \multicolumn{9}{c|}{$x_i \leq x_j$} \\
[.2cm]
$i$, $j$ & $-$ & 13, 14 & 18, 19 & 18, 19 & 9, 10 & $-$ & 6, 7 & 12, 13 & 18, 19 \\
\hline
\hline
Constraint & \multicolumn{9}{c|}{Tubes outer diameter is between the tubes inside diameter and the tubes spacing} \\
Indexation & $-$ & $-$ & $-$ & $c_{14},c_{15}$ & $c_{10},c_{11}$ & $-$ & $-$ & $-$ & $c_{16},c_{17}$ \\
[.2cm]
Equation & \multicolumn{9}{c|}{$x_i \leq x_j \leq x_k$} \\
[.2cm]
$i$, $j$, $k$ & $-$ & $-$ & $-$ & 22, 23, 20 & 13, 14, 11 & $-$ & $-$ & $-$ & 22, 23, 20 \\
\hline
\end{tabular}
\end{center}
\end{footnotesize}
\caption{Linear a~priori constraints in Instances~\solar{1.1} to~\solar{9.1} (\solar{10.1} is unconstrained).
A ``$-$'' means the constraint is not part of the instance.}
\label{tab-linear-apriori-constraints}
\end{table}

%
%
%
%
%
%
%
%

\subparagraph{Nonlinear a~priori constraints:}

Constraints details are shown in Table~\ref{tab-nonlinear-apriori-constraints}.
In \solar{2.1}, the left hand side of the heliostats field surface area constraint corresponds to the objective function to minimize (see Section~\ref{sec-obj}).

\begin{table}[H]
\renewcommand{\tabcolsep}{4.5pt}
\begin{footnotesize}
\begin{center}
\begin{tabular}{|l|ccccccccc|}
\hline
& \solar{1.1} & \solar{2.1} & \solar{3.1} & \solar{4.1} & \solar{5.1} & \solar{6.1} & \solar{7.1} & \solar{8.1} & \solar{9.1} \\
\hline
\hline
Constraint & \multicolumn{9}{c|} {Heliostats field surface area is lower than the maximum field surface area} \\
Indexation & $c_2$ & $c_1$ & $c_1$ & $c_1$ & $-$ & $-$ & $-$ & $c_1$ & $c_3$ \\
[.2cm]
Equation & \multicolumn{9}{c|}{$\pi x_3^2 (x_9^2-x_8^2)  x_7/180 \leq a$} \\
[.2cm]
$a$ [hectares] & 195 & 400 & 80 & 200 & $-$ & $-$ & $-$ & 400 & 500 \\
\hline
\hline
Constraint & \multicolumn{9}{c|}{Number of tubes in receiver fits inside receiver} \\
Indexation & $-$ & $c_{12}$ & $c_{11}$ & $c_{11}$ & $c_{7}$ & $-$ & $c_5$ & $c_7$ & $c_{12}$ \\
[.2cm]
Equation & \multicolumn{9}{c|}{$x_i x_j  \leq a\frac{\pi}{2}$} \\
[.2cm]
$i$, $j$, $a$ & $-$ & 11, 14, $x_5$ & 16, 19, $x_5$ & 16, 19, $x_5$ & 7, 10, 6 & $-$ & 4, 7, $x_2$ & 10, 13, $x_5$ & 16,19,$x_5$ \\
\hline
\end{tabular}
\end{center}
\end{footnotesize}
\caption{Nonlinear a~priori constraints in Instances~\solar{1.1} to~\solar{9.1} (\solar{10.1} is unconstrained).
A ``$-$'' means the constraint is not part of the instance.}
\label{tab-nonlinear-apriori-constraints}
\end{table}

%
%
%
%
%
%

\subparagraph{Simulation constraints:}

These constraints depend on some outputs of the simulation, described below.
Lengths, distances, tower size, tubes spacing and diameters are all in meters,
energy and parasitic load are expressed in~k$\W$h,
costs are in~\$,
surfaces in hectares,
temperatures in Kelvin,
and pressures in~$\mathrm{MPa}$.
The constraints are relaxable and quantifiable
(code QRSK in~\cite{LedWild2015}).
Moreover, some of these constraints correspond to stochastic outputs of the program.
This behavior is controlled with the {\tt -seed} option (see Section~\ref{subsec-stoch}).
Constraints details are shown in Table~\ref{tab-simulated-constraints}.

\begin{table}[htb!]
\renewcommand{\tabcolsep}{4.5pt}
\begin{footnotesize}
\begin{center}
\begin{tabular}{|l|ccccccccc|}
\hline
& \solar{1.1} & \solar{2.1} & \solar{3.1} & \solar{4.1} & \solar{5.1} & \solar{6.1} & \solar{7.1} & \solar{8.1} & \solar{9.1} \\
\hline
\hline
Constraint & \multicolumn{9}{c|}{The required number of heliostats ($x_6$) can fit in the heliostats field} \\
Indexation & $c_5$ & $c_6$ & $c_5$ & $c_5$ & $-$ & $-$ & $-$ & $c_4$ & $c_6$ \\
\hline
\hline
Constraint & \multicolumn{9}{c|}{The cost of plant is lower than the budget} \\
Indexation & $c_1$ & $c_3$ & $-$ & $-$ & $c_1$ & $-$ & $c_1$ & $-$ & $c_1$ \\
Budget limit & \$50M & \$300M & $-$ & $-$ & \$100M & $-$ & \$45M & $-$ & \$1.2B \\
\hline
\hline
Constraint & \multicolumn{9}{c|}{Compliance to demand} \\
Indexation & $-$ & $c_2$ & $c_2$ & $c_2$ & $-$ & $c_1$ & $-$ & $-$ & $-$ \\
Stochastic & $-$ & yes & yes & yes & $-$ & no & $-$ & $-$ & $-$ \\
\hline
\hline
Constraint & \multicolumn{9}{c|}{Minimal acceptable energy production} \\
Indexation & $-$ & $-$ & $-$ & $-$ & $-$ & $-$ & $-$ & $c_8$ & $c_2$ \\
Stochastic & $-$ & $-$ & $-$ & $-$ & $-$ & $-$ & $-$ & yes & yes \\
\hline
\hline
Constraint & \multicolumn{9}{c|}{Parasitic losses are lower than a ratio of the generated output} \\
Indexation & $-$ & $-$ & $-$ & $c_{13}$ & $c_9$ & $-$ & $c_6$ & $c_9$ & $c_{14}$ \\
Limit ratio & $-$ & $-$ & $-$ & 18\% & 18\% & $-$ & 3\% & 8\% & 20\% \\
Stochastic & $-$ & $-$ & $-$ & yes & no & $-$ & yes & yes & yes \\
\hline
\hline
Constraint & \multicolumn{9}{c|}{Storage is back at least at its original conditions} \\
Indexation & $-$ & $-$ & $c_{13}$ & $-$ & $-$ & $c_6$ & $-$ & $-$ & $-$ \\
Stochastic & $-$ & $-$ & yes & $-$ & $-$ & no & $-$ & $-$ & $-$ \\
\hline
\hline
Constraint & \multicolumn{9}{c|}{Molten salt temperature does not fall below the melting point in hot storage} \\
Indexation & $-$ & $c_8$ & $c_7$ & $c_7$ & $c_3$ & $c_3$ & $-$ & $-$ & $c_8$ \\
Stochastic & $-$ & yes & yes & yes & no & no & $-$ & $-$ & yes \\
\hline
\hline
Constraint & \multicolumn{9}{c|}{Molten salt temperature does not fall below the melting point in cold storage} \\
Indexation & $-$ & $c_9$ & $c_8$ & $c_8$ & $c_4$ & $c_4$ & $-$ & $-$ & $c_9$ \\
Stochastic & $-$ & yes & yes & yes & no & no & $-$ & $-$ & yes \\
\hline
\hline
Constraint & \multicolumn{9}{c|}{Molten salt temperature does not fall below the melting point in steam generator outlet} \\
Indexation & $-$ & $-$ & $c_9$ & $c_9$ & $c_5$ & $-$ & $-$ & $-$ & $c_{10}$ \\
Stochastic & $-$ & $-$ & no & yes & no & $-$ & $-$ & $-$ & yes \\
\hline
\hline
Constraint & \multicolumn{9}{c|}{Receiver outlet temperature is greater than the steam turbine inlet temperature} \\
Indexation & $-$ & $c_{12}$ & $c_{12}$ & $c_{12}$ & $c_8$ & $c_5$ & $-$ & $-$ & $c_{13}$ \\
\hline
\hline
Constraint & \multicolumn{9}{c|}{Pressure in receiver tubes is lower than yield pressure} \\
Indexation & $-$ & $c_7$ & $c_6$ & $c_6$ & $c_2$ & $c_2$ & $c_2$ & $c_5$ & $c_7$ \\
Stochastic & $-$ & yes & yes & yes & no & no & yes & yes & yes \\
\hline
\hline
Constraint & \multicolumn{9}{c|}{Pressure in steam generator tubes is lower than yield pressure} \\
Indexation & $-$ & $-$ & $-$ & $c_{16}$ & $c_{12}$ & $-$ & $-$ & $-$ & $c_{17}$ \\
\hline
\end{tabular}
\end{center}
\end{footnotesize}
\caption{Nonlinear a~priori constraints in Instances~\solar{1.1} to~\solar{9.1} (\solar{10.1} is unconstrained). A ``$-$'' means the constraint is not part of the instance. When stochasticity is not specified, constraints are deterministic.}
\label{tab-simulated-constraints}
\end{table}

\section[The SOLAR package for benchmarking]{The \solar{} package for benchmarking}
\label{sec-thepackage}

This section describes the \solar{} package version 1.0.
The correct way to refer to a specific instance requires the version number of the package.
For example, tests involving the third instance with the version 1.x of \solar{} should be denoted by \solar{3.1}: the first number is the instance, and the second is the main number of the version.


The complexity and nature of the \solar{} BBO problems depend on the selected instances and the options passed to the \solar{} framework. The number and nature of variables, the number of constraints and objectives have already been mentioned. In addition, for testing the agility of optimization methods to handle BBO problems, some features on the outputs can be present or not in \solar{} instances. These features are discussed in the following sections as such:
variables types in Section~\ref{subsec-var-types},
constraint types and failure to return valid outputs in Section~\ref{subsec-cstr-eval},
stochasticity in Section~\ref{subsec-stoch},
static surrogate models in Section~\ref{subsec-surrogates},
multifidelity in Section~\ref{subsec-multi-fid}.

\begin{figure}[ht!]
\begin{center}
\includegraphics[width=.9\textwidth, trim={0 8mm 0 2mm},clip]{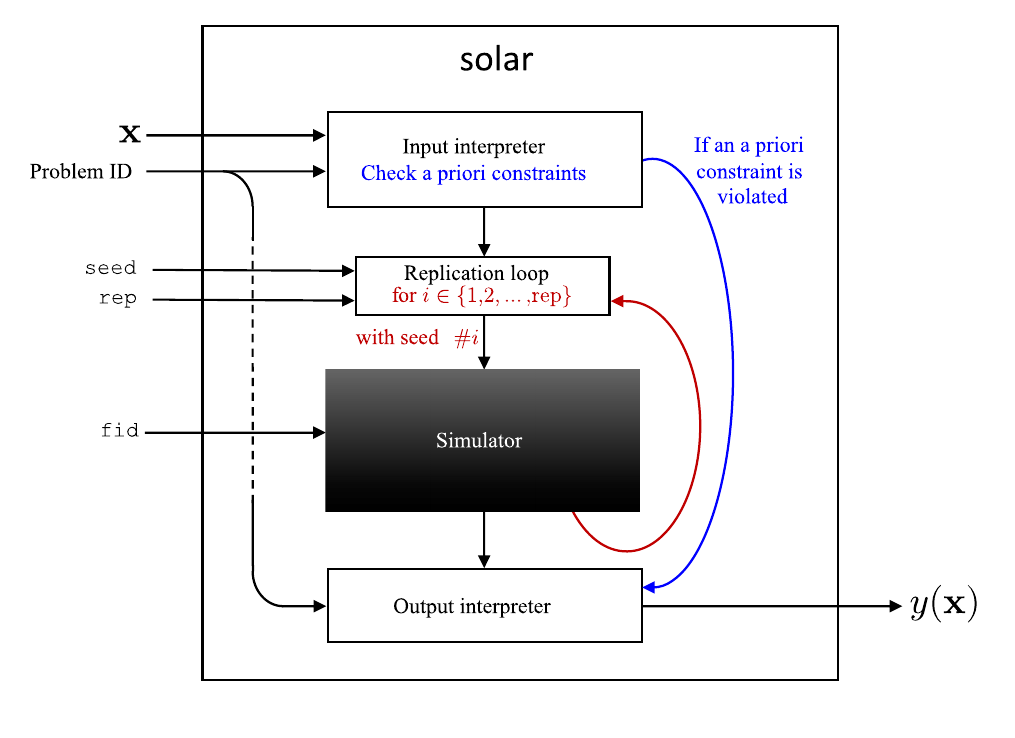}
\end{center}
\caption{The \solar{} framework. ``Problem~ID'' is the instance number between 1 and 10.}
\label{fig-framework}
\end{figure}

Figure~\ref{fig-framework} illustrates the inner working of the \solar{} framework.
For a given $\xx$, the first step consists in verifying if the a~priori constraints are satisfied 
     as detailed in Section~\ref{subsec-cstr-eval}.
 Then, the simulator is possibly called many times, with different random seeds 
     using the prescribed fidelity. In this case, the output interpreter returns averaged values of $F$ and of $C$ in the vector $y$.
In addition, the input file can contain several vectors $\xx$ and the output interpreter returns several vectors $y$. 

\subsection{Software utilization}
\label{subsec-software}

The~{\tt README} file from the GitHub package provides the necessary instructions to build \solar{} from source. As the software has been developed for comparable benchmarking, an in-house random number generator is used to guarantee that outputs are platform and compiler independent when calling the executable with the same settings. For this reason it is strongly encouraged to verify this behaviour by running the command ~{\tt solar -check} before conducting any benchmarking. Once the \solar{} version is stated to be valid\footnote{In the unlikely event that this is not the case, please contact the authors.}, all performance results can be compared to any studies which used the same \solar{} version for benchmarking.

Basic utilization information can be obtained from the executable by running the command ~{\tt solar~-h}. This describes how to run a simulation as well as the available options for stochasticity and fidelity and their default parameters; further information on these modes is given in their respective subsections. In addition, the best known values for single-objective instances (one replication, full fidelity, default random seed of zero) are reported.


Details on Instance~\verb|id| are obtained with the command \verb|solar -h id|. Table~\ref{tab-constraints-detailed} summarizes the information displayed by the command for all instances. The command also shows the suggested bounds on variables and a 
suggested starting point (feasible or infeasible) given in the \nomad parameters syntax. Additional information about the simulation can be obtained in verbose mode by adding \verb|-v| to the command. When the provided point satisfies all a~priori and hidden constraints, all the outputs are successfully computed by the simulator. In this case, the flag {\tt cnt\_eval=true} (``count evaluation''). This flag is displayed only in verbose mode.  

The \solar{} executable is obtained by building all the source code files. A more direct access to the \solar{} evaluator on a given instance is also possible (see the source code in the \verb|main_minimal| or \verb|main_minimal_2| functions in \verb|main.cpp|). The flag {\tt cnt\_eval} is an output of the evaluator call function \verb|eval_x|. 


\subsection{Mixed integer, categorical and unbound variables}
\label{subsec-var-types}

It is quite common for synthetic blackbox instances to have only real and bounded input variables. This however is only the case for Instances~\solar{6.1}
and~\solar{10.1}. Every other instance contains one or more integer inputs, as detailed in Tables~\ref{tab-vars1} and~\ref{tab-vars2}. If a non-integer value is provided in the input file for an integer variable, all the outputs returned by the code will be~$10^{20}$. Filtering or rounding values of input variables should be performed before launching \solar.   

Furthermore, Instances~\solar{3.1},~\solar{4.1},~\solar{5.1}, and~\solar{9.1} also contain categorical inputs. Whilst these take their values from a discrete set, they do not possess any ordering properties~\cite{AuDe01a}. A framework for modeling such variables in the context of BBO is presented in~\cite{G-2022-11}.

As an added difficulty, at least one integer input is always unbound, meaning approaches such as Latin hypercube Sampling which requires an upper and lower bound for every input variable cannot be directly applied and need to be adapted in some way.
Note that as described in Section~\ref{subsec-software}, calling the help command for a particular instance provides, among others, a suggested starting point for optimization. The suggested starting value for unbound input variables can be used as a guide for what input range to explore, although it provides no guarantees in terms of proximity to an optimal point.

\begin{table}[ht!]
\renewcommand{\tabcolsep}{4.5pt}
\begin{footnotesize}
\begin{center}
\begin{tabular}{|l|ccc|}
\hline
Instance & A priori & Simulated deterministic & Simulated stochastic \\
\hline
\hline
\solar{1.1}  & $c_2,c_3,c_4$                             & $c_1,c_5$                              & $f_1$                             \\
\solar{2.1}  & $f_1,c_1,c_4,c_5,c_{10},c_{11}$             & $c_3,c_6,c_{12}$                       & $c_2,c_7,c_8,c_9$               \\
\solar{3.1}  & $c_1,c_3,c_4,c_{10},c_{11}$               & $f_1,c_5,c_9,c_{12}$                     & $c_2,c_6,c_7,c_8,c_{13}$        \\
\solar{4.1}  & $c_1,c_3,c_4,c_{10},c_{11},c_{14},c_{15}$ & $f_1,c_5,c_{12},c_{16}$                  & $c_2,c_6,c_7,c_8,c_9,c_{13}$    \\
\solar{5.1}  & $c_6,c_7,c_{10},c_{11}$                   & $f_1,c_1,c_2,c_3,c_4,c_5,c_8,c_9,c_{12}$ &                                 \\
\solar{6.1}  &                                           & $f_1,c_1,c_2,c_3,c_4,c_5,c_6$            &                                 \\
\solar{7.1}  & $c_3,c_5$                                 & $c_1,c_4$                              & $f_1,c_2,c_6$                     \\
\solar{8.1}  & $c_1,c_2,c_3,c_6,c_7$                     & $f_2,c_4$                              & $f_1,c_5,c_8,c_9$               \\
\solar{9.1}  & $c_3,c_4,c_5,c_{11},c_{12},c_{15},c_{16}$ & $c_1,c_6,c_{13},c_{17}$                & $f_1,f_2,c_2,c_7,c_8,c_9,c_{10},c_{14}$ \\
\solar{10.1} &                                           & $f_1$                                    &                                 \\
\hline
\end{tabular}
\end{center}
\end{footnotesize}
\caption{Classification of the outputs of the ten blackbox instances as either a~priori, simulated and deterministic or simulated and stochastic.}
\label{tab-constraints-detailed}
\end{table}

\subsection{Constraints evaluation}
\label{subsec-cstr-eval}

The \solar{} package includes a~priori (QUAK)~\cite{LedWild2015} constraints with given expressions (see Section~\ref{sec-cstr}). The specific outputs are shown in Table~\ref{tab-constraints-detailed}. They are fast to evaluate compared to simulation constraints. Hence, the simulator does not start if at least one is violated. In this case, the values for the a~priori constraints are returned but those for the remaining simulation constraints are set to~$10^{20}$. This evaluation should not be counted during an optimization, as indicated by the raised flag {\tt cnt\_eval=false} (available only in verbose mode or when calling directly the \solar{} evaluator). Table~\ref{tab-constraints-feasibility} gives an estimation of how frequently different types of constraints are violated for each instance.

\begin{table}[ht!]
\renewcommand{\tabcolsep}{4.5pt}
\begin{footnotesize}
\begin{center}
\begin{tabular}{|l|r|rrr|rrr|}
\hline
\multirow{2}{*}{Instance} & Eval. (nb) & \multicolumn{3}{c|}{LH sampling (10,000 points)} & \multicolumn{3}{c|}{\nomad 3} \\
& $=400n$ & feas. AP (\%) & feas. (\%) & hidden (\%) & feas. AP (\%) & feas. (\%) & hidden (\%) \\
\hline
\hline
\solar{1.1}  &  3,600 & 29.0 & 4.28 & 0 & 96.9 & 73.0 & 0.417  \\
\solar{2.1}  &  5,600 & 17.9 & 0 & 0 & 91.8 & 6.89 & 0  \\
\solar{3.1}  &  8,000 & 8.63 & 0 & 0 & 94.9 & 11.4 & 0.100  \\
\solar{4.1}  & 11,600 & 0.112 & 0 & 0 & 84.8 & 28.0 & 0.0603 \\
\solar{5.1}  &  3,000 & 1.60 & 0.10 & 0 & 88.2 & 14.3 & 0 \\
\solar{6.1}  &  2,000 & 90.1 & 4.95 & 0 & 100 & 40.0 & 0.557 \\
\solar{7.1}  &  2,800 & 46.4 & 14.2 & 0 & 82.7 & 50.9 & 0.143 \\
\solar{8.1}  &  5,200 & 0.942 & 0 & 0 & 80.2 & 60.3 & 0.135 \\
\solar{9.1}  & 11,600 & 0.845 & 0 & 0.0948 & 64.6 & 9.82 & 0.793 \\
\solar{10.1} &  2,000 & 100 & 100 & 0 & 100 & 100 & 0 \\
\hline
\end{tabular}
\end{center}
\end{footnotesize}
\caption{Percentage of feasibility of a~priori constraints (``feas. AP''), of all constraints (``feas.'') and of hidden constraints (``hidden'') for Latin hypercube sampling and for \nomad 3~\cite{Le09b}. For comparability between instances, an evaluation budget of~400 times the instance dimension is given to each instance. The \nomad optimization of \solar{6.1} converged and stopped after 1,976 evaluations, and~150$n$ evaluations are given to \solar{5.1} because its evaluation time is high.}
\label{tab-constraints-feasibility}
\end{table}

Hidden (NUSH)~\cite{LedWild2015} constraints are present in some instances. Some bugs may still be present in the code, as is often the case in real blackbox problems. Moreover, some a~priori constraints may not be explicitly given. For example, when $x_3 x_8 \neq 0$ is not verified in \solar{1.1}, a hidden constraint is violated. Some hidden constraints can even stop an evaluation before any simulation is launched (see Figure~\ref{fig-framework}), causing the flag {\tt cnt\_eval=false} to be raised. When a hidden constraint is violated, the \solar{} evaluation does not complete, and some simulated outputs return the value~$10^{20}$. Not all outputs are necessarily affected by a violated hidden constraint. Note that for most cases, the closer to 1 the fidelity is, the less likely it is to violate a hidden constraint. This is not always the case, as for \solar{3.1}, it is possible to find a point such that a hidden constraint is violated only with fidelity 1. To benchmark a blackbox with hidden constraints, the authors recommend \solar{3.1} with lower fidelities. This way, violating a hidden constraint is still relatively rare. To benchmark in a situation of highly occurring hidden constraints violation, the authors recommend the STYRENE problem~\cite{AuBeLe08}\footnote{Available at  \href{https://github.com/bbopt/styrene}{\url{https://github.com/bbopt/styrene}}}.

\subsection{Stochasticity}
\label{subsec-stoch}

Table~\ref{tab-constraints-detailed} shows that seven of the ten instances contain at least one stochastic output. The collection of instances of this type provides a variety of benchmarks for algorithm development with different characteristics. Specifically, \solar{1.1} provides a classical stochastic optimization benchmark where only the objective function is stochastic, but some deterministic constraints also exist. Instance~\solar{7.1} contains a stochastic objective function as well as some stochastic constraints.
Instances~\solar{2.1},~\solar{3.1} and~\solar{4.1}, on the other hand, provide less traditional benchmarks where the objective function is deterministic, but some of the constraints are stochastic. Finally, for an added layer of benchmark complexity, Instances~\solar{8.1} and~\solar{9.1} are defined
with two objectives and stochastic outputs.

\subsubsection{The {\tt seed} argument (random seed)}

\begin{figure}[ht!]
    \begin{center}
        \includegraphics[width=1\textwidth]{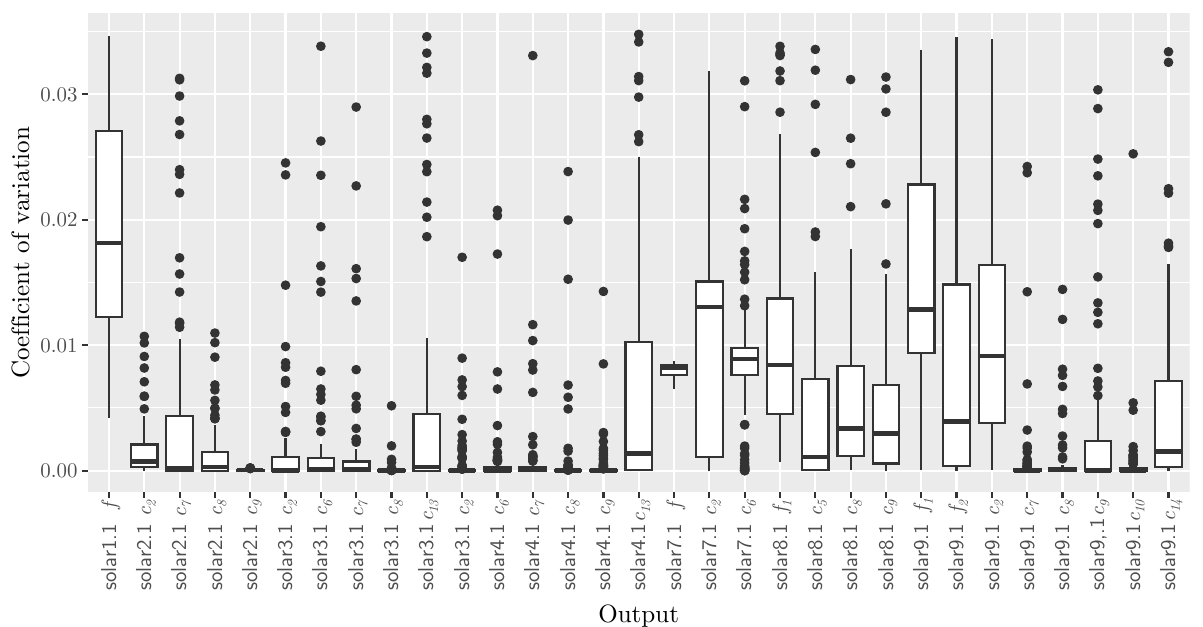}
    \end{center}
    \caption{Distribution of the coefficient of variation of the stochastic outputs of all instances. Each box plot shows the coefficient of variation calculated at 100 randomly selected points, calculated for 100 seeds each. Each box plot shows the median (thick black line), the interquartile range (rectangle), and the min and max range (vertical lines) excluding outliers (black points). Note that a few outliers outside of the plot range are not shown here for clarity.}
    \label{fig-stochastic}
\end{figure}

For reproducibility purposes, every stochastic behaviour stems from a seed, controlled with the~{\tt -seed} argument passed to the \solar{} command. The custom pseudo-random number generator used by \solar{} also guarantees stochastic reproducibility for the same inputs and specified seed. This behaviour is verified when running the command~{\tt solar -check} after compiling \solar. The default seed value is 0, and thus no stochastic behaviour will be observed between executions if this parameter is not specifically modified. To enable stochasticity, set {\tt -seed=diff} to randomly generate a different seed each time, leading to non-reproducible behaviour. Users can also choose to seed the random generator with a chosen integer {\tt s} with the argument {\tt -seed=s}. Note that another potential use of fixing the seed is the creation of different deterministic instances of the same simulation. Figure~\ref{fig-stochastic} shows how stochasticity leads to different levels of variability for the stochastic outputs of all instances. Even within the same instance, variability can be very different among the outputs; in \solar{9.1} for example Constraints~7,~8 and~10 exhibit little variations compared to Constraint~2. Choosing how to approach different stochastic levels for each of the outputs simultaneously is up to the technique being assessed.

\subsubsection{The {\tt rep} argument (replications)}

A common algorithmic approach to handle stochastic outputs involves repeatedly sampling the same location and taking the output average. In order to aid practitioners in the use of \solar{}, this approach can be directly queried from the simulator via the argument~{\tt -rep=r} with~${\tt r} \in\N^+$ passed to the \solar{} command. Passing this parameter will run the specified simulation a total of~{\tt r} times and return only the average of the outputs. The default~{\tt -rep} value is 1. Additionally,  when the number of repetitions is greater than one, the~{\tt -seed} argument is now used to generate different repetitions. This means that running a \solar{} instance with~{\tt -seed=s} and~{\tt -rep=r} will always result in the same repetitions, the same output values. Again, this is true across all platforms and compilers.

Alternatively, another algorithmic approach is to repeatedly sample a point until
a stopping rule is satisfied.
The Gauss-distribution criteria, described in~\cite{BiWaBrPo22}, has been implemented.
In practice, it is activated by choosing ${\tt r}$ as a real in $(0,1)$ for the
{\tt -rep=r} argument. In this mode, ${\tt r}$ represents the probability that all the
stochastic outputs are stabilized within 3 significant digits. A replication number limit of~$5\times10^4$ ensures the existence of an upper bound on computational time.
Table~\ref{tab-real-r} shows that closer to an optimum, the number of replications needed to reach each probability is generally considerably higher, although exceptions are present. Notice that even with probabilities as low as~$0.05$, it is possible to reach the replication limit. Finally one must be careful when using this parameter within a time budget. Notice that \solar{2.1} at the best known point with probability~$0.5$ took almost 12 days to reach~$5\times10^4$ replications.

\begin{table}[ht!]
\renewcommand{\tabcolsep}{4.5pt}
\begin{footnotesize}
\begin{center}
\begin{tabular}{|lc|rrrrr|rrrrr|}
\hline
\multirow{2}{*}{Instance} & \multirow{2}{*}{Input} & \multicolumn{5}{c|}{nb of replications for multiple {\tt r} values} & \multicolumn{5}{c|}{CPU time (seconds) for multiple {\tt r} values} \\
 &  & 0.05 & 0.25 & 0.5 & 0.75 & 0.95 & 0.05 & 0.25 & 0.5 & 0.75 & 0.95 \\
\hline
\hline
\multirow{2}{*}{\solar{1.1}} & $\xx^0$ & 7 & 125 & 558 & 1,558 & 4,299 & 0 & 9 & 58 & 275 & 1,567 \\
                     & $\xx^{\star}$ & 10 & 36 & 94 & 232 & 678 & 84 & 324 & 825 & 2,076 & 5,824 \\
\multirow{2}{*}{\solar{2.1}} & $\xx^0$ & 15 & 164 & 688 & 2,033 & 5,837 & 136 & 1,497 & 6,480 & 19,574 & 58,875 \\
                     & $\xx^{\star}$ & 884 & 22,571 & 50k & 50k & 50k & 14,059 & 468,595 & 1,017k & 1,018k & 1,004k \\
\multirow{2}{*}{\solar{3.1}} & $\xx^0$ & 7 & 29 & 194 & 562 & 1,582 & 12 & 50 & 324 & 1,025 & 2,996 \\
                     & $\xx^{\star}$ & 13,153 & 50k & 50k & 50k & 50k & 38,904 & 147,781 & 142,858 & $-$ & $-$ \\
\multirow{2}{*}{\solar{4.1}} & $\xx^0$ & 7 & 7 & 7 & 7 & 7 & 13 & 13 & 13 & 13 & 13 \\
                     & $\xx^{\star}$ & 3,723 & 50k & 50k & 50k & 50k & 13,578 & 266,764 & $-$ & $-$ & 242,440 \\
\multirow{2}{*}{\solar{7.1}} & $\xx^0$ & 7 & 16 & 60 & 148 & 369 & 25 & 57 & 212 & 525 & 1,311 \\
                     & $\xx^{\star}$ & 50k & 50k & 50k & 50k & 50k & 256,757 & $-$ & $-$ & $-$ & $-$ \\
\multirow{2}{*}{\solar{8.1}} & $\xx^0$ & 7 & 7 & 38 & 122 & 415 & 43 & 43 & 232 & 789 & 2,560 \\
                     & $\xx^a$ & 7 & 7 & 12 & 18 & 65 & 49 & 49 & 87 & 128 & 450 \\
\multirow{2}{*}{\solar{9.1}} & $\xx^0$ & 26 & 285 & 1,164 & 3,366 & 9,591 & 70 & 774 & 3,237 & 9,960 & 33,161 \\
                     & $\xx^b$ & 8 & 85 & 486 & 1,344 & 3,681 & 31 & 369 & 2,054 & 6,407 & 21,180 \\
\hline
\end{tabular}
\end{center}
\end{footnotesize}
\caption{Number of replications and CPU time (with an Apple M2 max) for some instances when using five different values of~{\tt -rep=r} in~$(0,1)$. Input~$\xx^{\star}$ denotes the best known point at the time of writing, whereas~$\xx^0$ denotes the starting point suggested for the problem,
far from~$\xx^{\star}$ (these points are further described in Table~\ref{tab-x0xs}). For multiobjective instances, since no known point dominates all the others, random points~$\xx^a$ and~$\xx^b$ are shown. For a probability parameter higher than another where the replication number limit was reached at the same point, the simulation is sometimes not run since it would surely reach the limit as well, and a~``$-$'' is shown.}
\label{tab-real-r}
\end{table}


\subsection{Surrogate models}
\label{subsec-surrogates}


In the present context of Problem~\eqref{eq-pb}, a surrogate is a simplified version of the functions $F$ and $C$ that are replaced by approximations which are less computationally expensive.
These surrogates are constructed by relaxing some internal parameters within the simulation.
For example, a surrogate may trace fewer rays, or sample the performance of the plant at fewer instants over the time window than the real blackbox.


Some solvers take advantage of such alternative models to conduct the optimization, 
	and use the original model only to verify promising solutions found  with the surrogate.
The surrogate management framework~\cite{BoDeFrSeToTr99a}
 proposes ways to exploit  surrogates in optimization  problems.


A surrogate model used by an optimization algorithm can be dynamic, meaning the algorithm builds its own model and adjusts it during the optimization, or static, meaning it is given to the algorithm beforehand. With the \solar{} package, there are different ways to create static surrogate models. First, by lowering the fidelity, a wide range of static surrogates with variable precision are created. Second, by considering a chosen number of replications for a stochastic blackbox to be the truth, any lower number of replications results in a static surrogate. A static surrogate model can also be achieved by lowering both the fidelity and the number of replications.

\subsection{Multifidelity}
\label{subsec-multi-fid}

\begin{figure}[ht!]
    \begin{center}
        \includegraphics[width=0.7\textwidth]{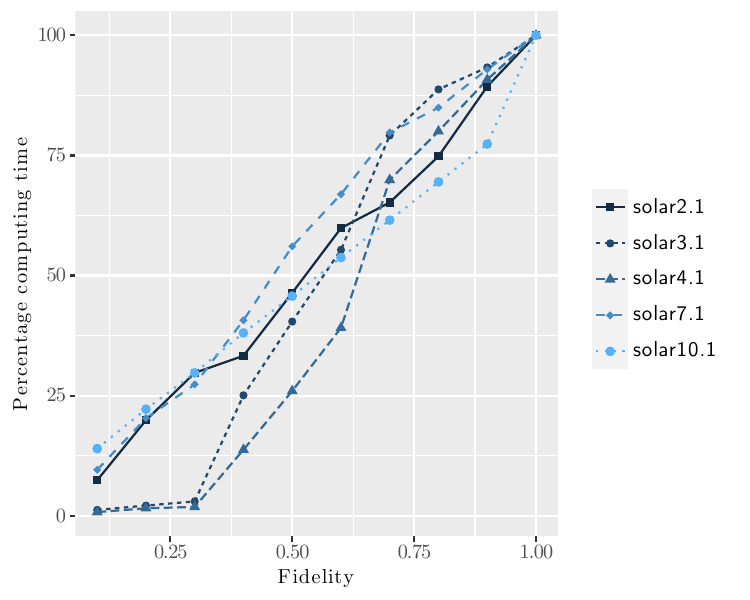}
    \end{center}
    \caption{Computation time for different fidelities on
    Instances~\solar{2.1},~\solar{3.1},~\solar{4.1},~\solar{7.1}, and~\solar{10.1},
    with respect to a fidelity of 1. 
    }
    \label{fig-multiFidelityTimes}
\end{figure}

With Instances~\solar{2.1},~\solar{3.1},~\solar{4.1},~\solar{7.1},~\solar{8.1},~\solar{9.1}, and~\solar{10.1}, a fidelity argument~{\tt -fid=$\phi$} can be passed to the \solar{} command, with~$\phi\in[0,1]$. The default~$\phi$ value is 1. 
Lowering the fidelity generally leads to faster computational times, at the cost of degraded accuracy of the simulation results. As \solar{} is a simulation engine, both the time and accuracy reduction do not follow a specific formula. With~{\tt -fid=0}, only the a~priori outputs are evaluated, while the others are set to $10^{20}$. Figure~\ref{fig-multiFidelityTimes} shows the relationship between fidelity and execution time on the single-objective multifidelity instances.

To focus on simulation time, only points that satisfy all a~priori constraints are considered. As it can be seen, the percentage reduction in time is roughly linear with the reduction in fidelity for all instances.


In all constrained instances where the fidelity parameter is available, a number of constraints are affected by the~{\tt -fid} parameter. This is a unique feature among BBO benchmarking problems, despite the fact that many industrial BBO problems have this property~\cite{AlAuDiLedLe23}.
Trivially, a set of bi-fidelity benchmarks similar to~\cite{To2015} can easily be generated by selecting pairs of fidelities.

\begin{figure}[ht!]
    \begin{center}
        \includegraphics[width=0.8\textwidth]{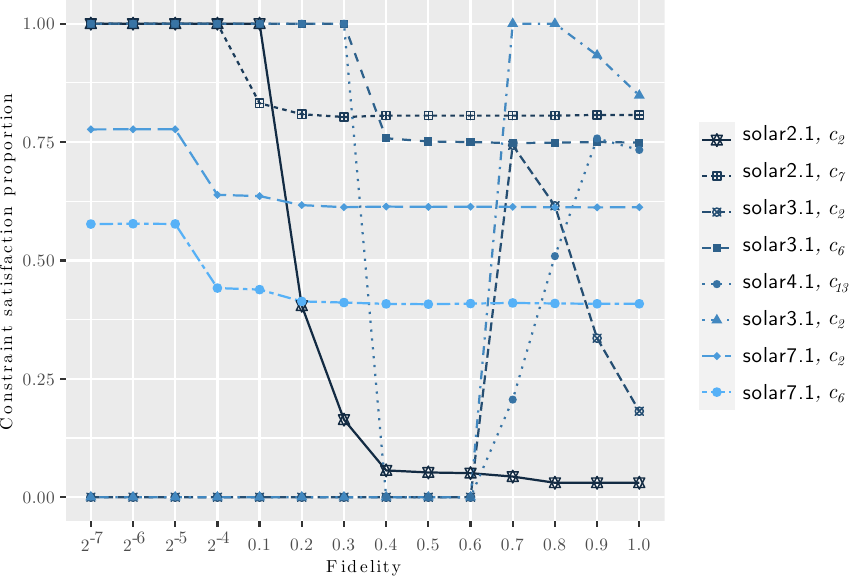}
    \end{center}
    \caption{Proportion of points within a small radius of the best known feasible solutions which satisfy particular constraints at different fidelities. 
    }
    \label{fig-multiFidelityConstraints}
\end{figure}

Figure~\ref{fig-multiFidelityConstraints} shows how lowering the fidelity can lead to the incorrect labeling of whether a constraint is violated or not. Certain constraints such as constraint 6 for \solar{7.1} are not much impacted, whereas others such as constraint 13 for \solar{4.1} is only reliable for a fidelity of 0.9 or higher.

\begin{figure}[ht!]
    \begin{center}
        \includegraphics[width=0.6\textwidth]{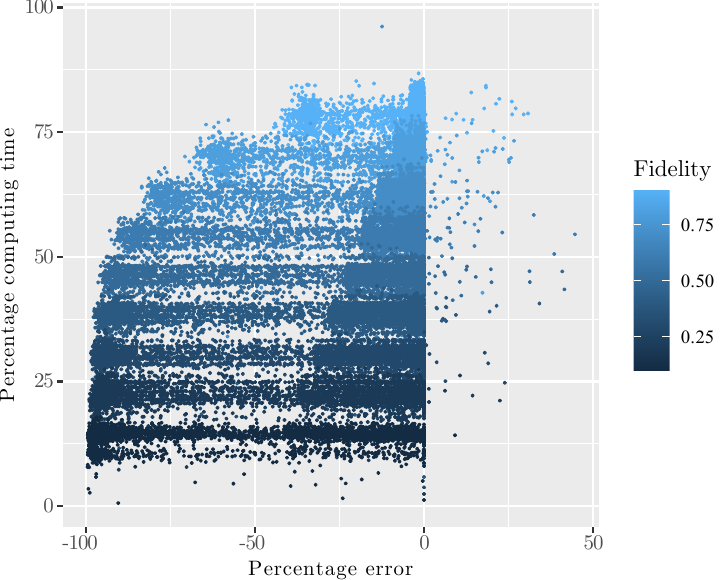}
    \end{center}
    \caption{Computational time and error on \solar{10.1} objective value on sampled points using different fidelities. Each sampled point is associated with a colored dot on the plot at coordinate $(x,y)$. The $y$-value denotes the percentage time taken compared to running the same sample point with fidelity of~1, and the $x$-value denotes the error relative to the objective function value of the same sample point with fidelity of~1. Note that a negative error denotes an underestimation, and positive error a overestimation of the objective function.}
    \label{fig-multiFidelityErrorAndTime}
\end{figure}

Figure~\ref{fig-multiFidelityErrorAndTime} showcases on \solar{10.1} the impact fidelity has both on computational time and percentage error compared to using a fidelity of~1. Instance~\solar{10.1} is an unconstrained problem, for which modifying the fidelity impacts only the objective function. As the figure shows, whilst providing no guarantees, increasing the fidelity leads to a consistent increase in time as evidenced by points of the same fidelity lying on the same height. At the same time, an increase in fidelity leads to a consistent decrease in error, as evidenced by a reduction in the width in which points of higher fidelities lie. Interestingly, the vast majority of points feature an underestimation of the objective value.

\subsection{Examples of optimization results}

Part of the challenge of building a new blackbox for BBO benchmarking is making sure that the resulting optimization problems meet some basic requirements: they need to have non-trivial optimal solutions, and they need to be acceptably functional. For instance, if a problem would be to minimize the cost of the designed power plant, it has to contain constraints sufficient to ensure that an optimization algorithm would not cling to a trivial solution such as a power plant with no heliostats at all that has a very low cost but also produces absolutely no energy. 

The best known points, at the time of writing, for the single-objective instances, are noted $\xx^{\star}$ and have been collected in the recent years with the help of several beta-testers (they are named in the acknowledgments).
Along with the suggested starting points, noted $\xx^0$, they are summarized in Table~\ref{tab-x0xs}.
The best known values $f_1(\xx^{\star})$ are also indicated in the GitHub repository, but not the coordinates $\xx^{\star}$ themselves.

In all cases, the best known solutions are non-trivial, as noticeable by the low number of active bounds at $\xx^{\star}$ and the number of saturated inequality constraints. Other constraints closed to saturation are not reported in the table but this may be an indication that the objective values may be further reduced although probably not easily.
Moreover, in Instances~\solar{1.1},~\solar{3.1} and~\solar{7.1}, the objective value at $\xx^{\star}$ has been significantly reduced compared with $f_1(\xx^0)$.

\begin{table}[H]
\renewcommand{\tabcolsep}{3.75pt}
\begin{footnotesize}
\begin{center}
\begin{tabular}{|r|rrrrrrrr|}
\hline
& \solar{1.1} & \solar{2.1} & \solar{3.1} & \solar{4.1} &
\solar{5.1} & \solar{6.1} & \solar{7.1} & \solar{10.1} \\
\hline
\hline
 $f_1(\xx^0)$  & -122,506  & 753,527 & 107,541,652 & 78,622,308 &  -30.0615 & 4,136,232 &-2,768 & 1,880\\
 $f_1(\xx^{\star}$) & -902,504 & 841,840 & 70,813,885 & 108,197,236 &  -28.8817 & 43,954,935 & -4,973 & 42 \\
 \hline
 \hline
  viol. constr. at $\xx^0$ & feasible & $c_2$, $c_6$, $c_7$ & $c_2$ & $c_2$, $c_5$          &
   $c_2$,   $c_5$,  $c_9$ & $c_1$, $c_6$ & $c_2$, $c_6$ & feasible \\
active constr. at $\xx^0$ & $c_5$ & $-$                   & $c_5$ & $-$                     &
 $c_4$ & $-$ & $-$ & $-$ \\
 active bounds at $\xx^0$ & $-$     & $u_7$               & $-$     & $u_{21}$, $\ell_{27}$ &
 $\ell_{16}$, $\ell_{20}$ & $-$ & $-$ & $-$ \\
\hline
\hline
active constr. at $\xx^{\star}$ & $c_5$ & $c_2$, $c_6$ & $c_2$, $c_5$, $c_8$ & $c_2$, $c_5$ & $-$ & $c_1$ & $c_4$ & $-$ \\
active bounds at $\xx^{\star}$  & $-$     & $-$            & $\ell_{14}$         & $u_{15}$     & $u_6$, $\ell_{16}$, $\ell_{19}$ & $u_1$, $u_3$, $\ell_5$ & - & $\ell_5$ \\
\hline
\end{tabular}
\end{center}
\end{footnotesize}
\caption{Comparison of the values of the initial points ($\xx^0$) with the values of the best known points ($\xx^{\star}$), for the single-objective instances. These values are rounded, and when negative, they correspond to a maximization problem.
The lines ``viol. constr.'', ``active constr.'', and ``active bounds''
indicate which constraints or bounds are active or violated, for $\xx^0$ or $\xx^{\star}$.}
\label{tab-x0xs}
\end{table}



Data profiles obtained on \solar{6.1} with 30 different starting points generated by Latin Hypercube sampling show that the optimization solvers tested (NOMAD 3~\cite{Le09b}, NOMAD 4~\cite{nomad4paper} and CMA-ES~\cite{cmaes}) do not systematically produce the best known solution.\footnote{Data available at the \href{https://github.com/bbopt/solar/tree/master/tests/6_MINCOST_TS}{\solar{} GitHub repository}}  

\begin{figure}[ht!]
    \begin{center}
        \includegraphics[width=1.0\textwidth]{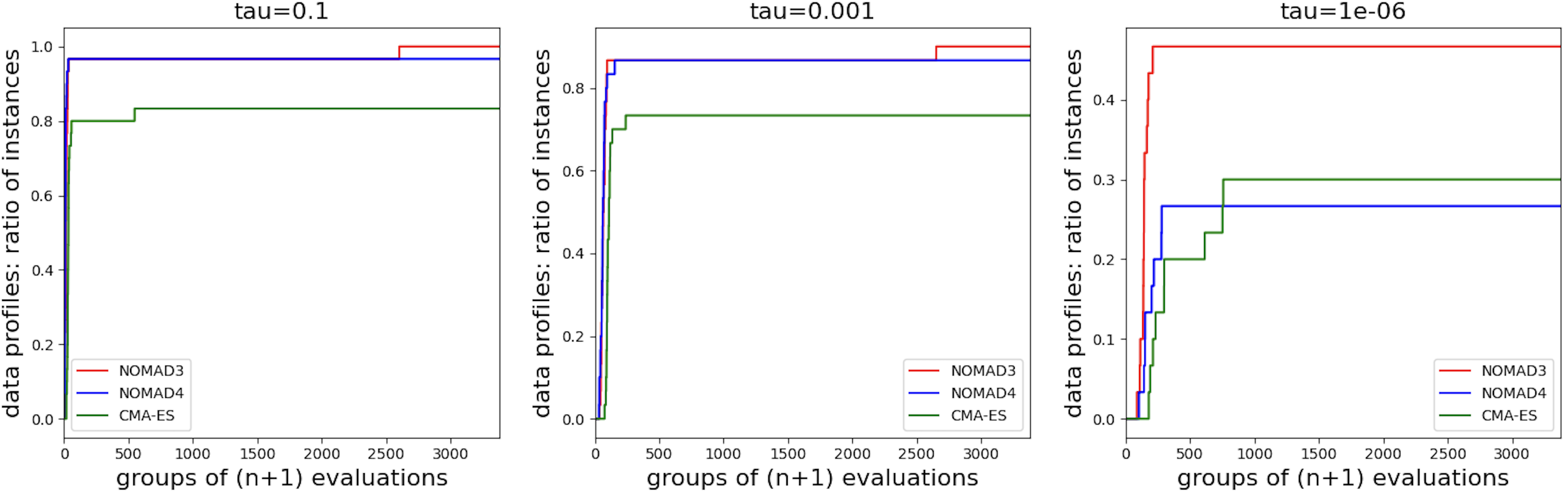}
    \end{center}
    \caption{Data profiles obtained on \solar{6.1} using 30 different starting points generated by Latin Hypercube sampling with the NOMAD~3~\cite{Le09b}, NOMAD~4~\cite{nomad4paper} and CMA-ES~\cite{cmaes} solvers.}
    \label{fig-dp-solar6}
\end{figure}

\begin{figure}[ht!]
    \begin{center}
        \includegraphics[width=0.9\textwidth]{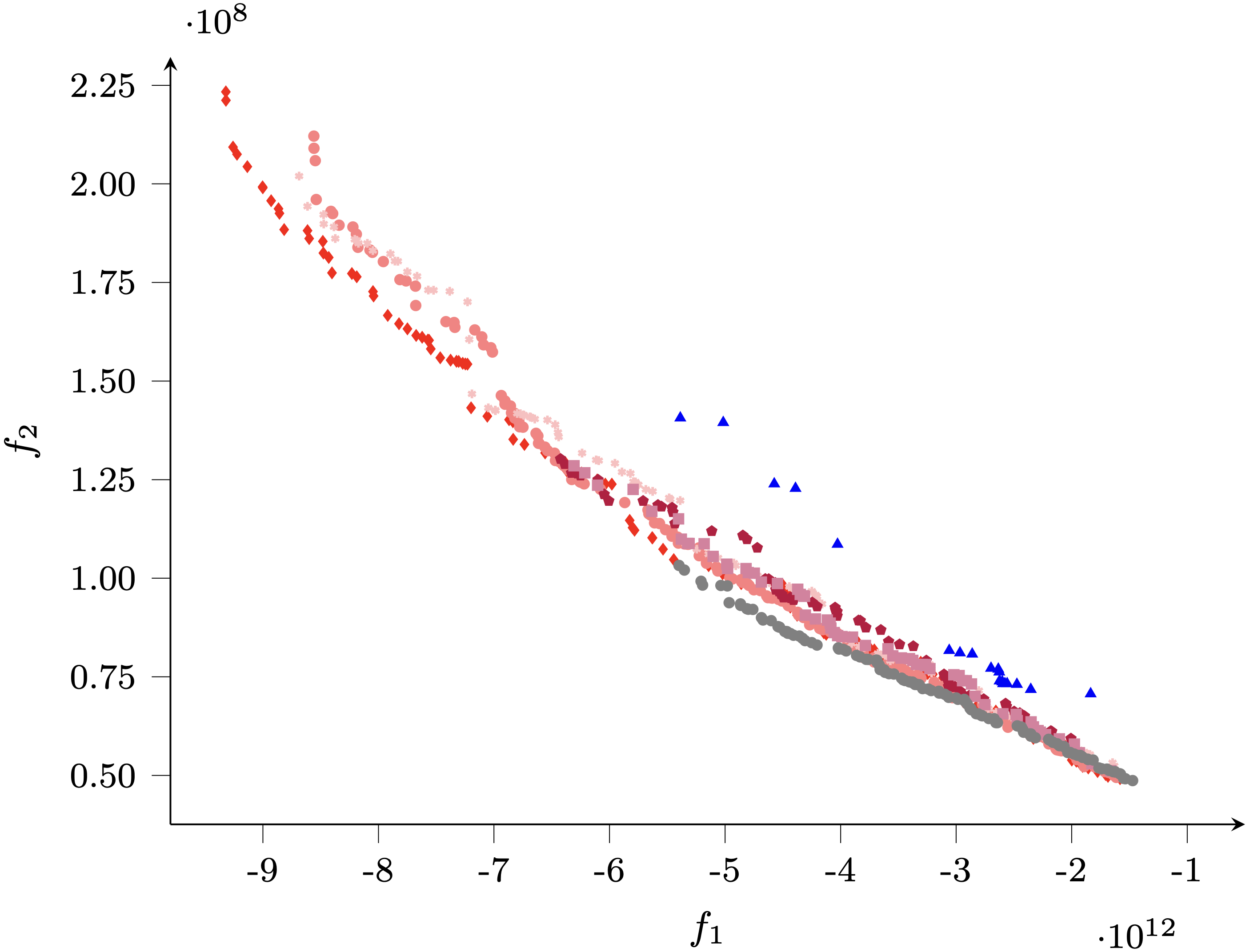}
    \end{center}
    \caption{Examples of non-dominated solutions for \solar{8.1} obtained during multiobjective optimizations conducted in~\cite{G-2022-10} using various solvers and a maximal budget of 5,000 evaluations.}
    \label{fig-multiObjective}
\end{figure}

As shown in Table~\ref{tab-constraints-detailed}, \solar{} provides two
multiobjective instances for benchmarking, namely~\solar{8.1} and~\solar{9.1}.
Figure~\ref{fig-multiObjective} displays examples of non-dominated solutions of both instances when using the default seed, fidelity and replication inputs. These plots give an indication of the characteristics of the respective Pareto fronts. The non-dominated solutions obtained by different solvers are distinct, suggesting that \solar{8.1} and \solar{9.1} are not trivial problems.

\subsection{Sensitivity of optimization solutions}

The best solutions obtained during the preliminary tests shown previously often feature saturated or close to saturation inequality constraints. Navigating the design space in a region around best solution presents algorithmic and numerical challenges for most solvers. The sensitivity of constraints feasibility around such points can be a good indicator of the difficulty of the task.  

\begin{figure}[ht!]
    \begin{center}
        \includegraphics[width=0.9\textwidth]{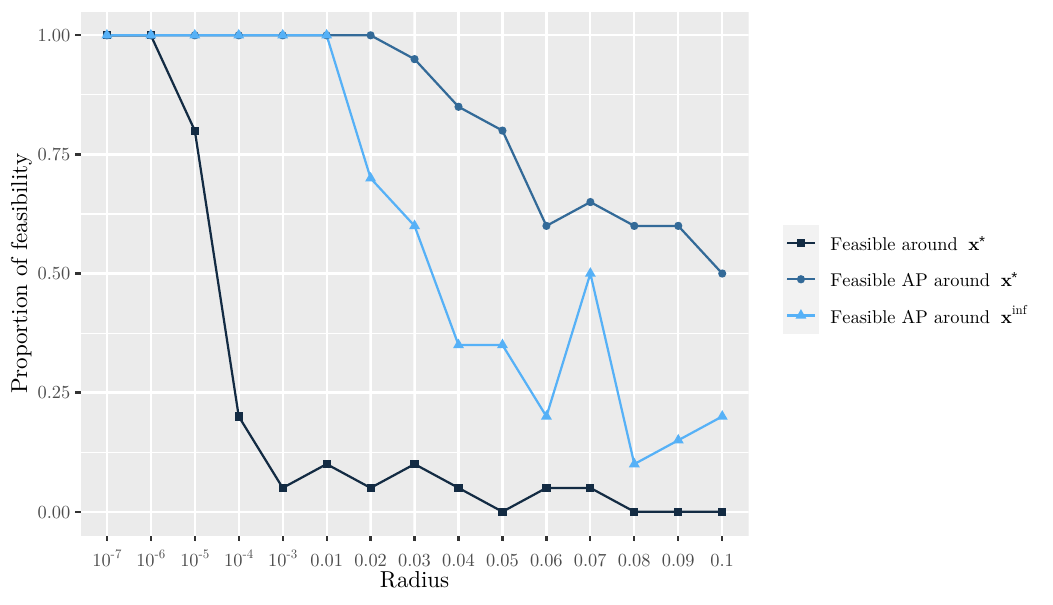}
    \end{center}
    \caption{Plot showing how the feasibility of \solar{4.1} points is sensitive to small changes. For different radii~$r$, ``Feasible around~$\xx^{\star}$'' shows the proportion of feasible points sampled in~$B_r(\xx^{\star})$, ``Feasible AP around~$\xx^{\star}$'' shows the proportion of points sampled in~$B_r(\xx^{\star})$ that satisfy all a~priori constraints, and ``Feasible AP around~$\xx^{inf}$'' shows the proportion of points sampled in~$B_r(\xx^{inf})$ that satisfy all a~priori constraints. Here,~$\xx^{\star}$ is the best known solution for \solar{4.1} at the time of writing and~$\xx^{inf}$ is an infeasible point.}
    \label{fig-sensibilityFeasibility}
\end{figure}

Figure~\ref{fig-sensibilityFeasibility} shows the proportion of feasible points (all constraints or only a~priori constraints) around two types of points obtained with \solar{4.1}. When considering all constraints, it can be noted that achieving feasibility requires to be very close to $\xx^{\star}$.

\begin{figure}[ht!]
    \begin{center}
        \includegraphics[width=0.7\textwidth]{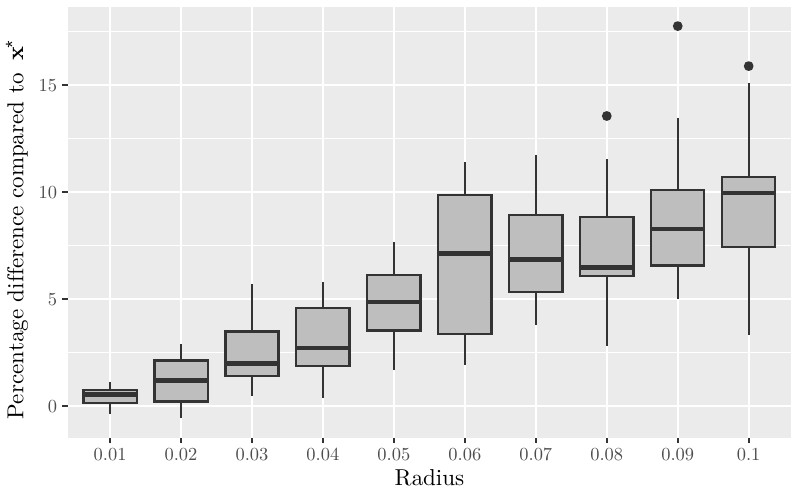}
    \end{center}
    \caption{Percentage change in the objective function of \solar{10.1} at different distances from a best known solution $\xx^{\star}$.}
    \label{fig-sensigilityPoint}
\end{figure}

In a realistic optimization problem, the direction and amplitude of change in objective function value play a role in a solver progress.
If objective function reduction is possible close to current best solution 
maybe there is room for improvement while still maintaining saturated inequality constraints. 
When solver stops, for an optimization solver end user, this can be an incentive to relax saturated inequality constraints (if possible) to obtain more gain.
Figure~\ref{fig-sensigilityPoint} shows the sensitivity of the objective function of \solar{10.1} near $\xx^{\star}$. Hence, we believe that \solar{10.1} is non trivial with characteristics similar to what is found in real life optimization problems and represents a good challenge for solvers. 

\section{Discussion}

There are very few industrial or real-life publicly available BBO problems.
This works addresses this issue by introducing \solar, a collection of benchmark problem instances for BBO.
\solar{} contains ten instances of a blackbox that simulates a concentrated solar power plant.
The instances differ in terms of variables, constraints, stochasticity and number of objective functions.
\solar{} is coded in standard {\sf C++} and is publicly available as a stand-alone open source package on GitHub.
Practitioners and algorithm developers are encouraged to test their methods with this blackbox and
to participate by sending their best optimization logs and results to
\href{mailto:nomad@gerad.ca?subject=Results\%20for\%20the\%20Solar\%20benchmark}{\tt nomad@gerad.ca}
in order to update the database of results on the \solar{} GitHub page.
Future work includes the definition and release of new instances of the blackbox to the community.

%

\section*{Acknowledgments}

The authors would like to thank the following persons for having tested the preliminaries versions of \solar,
and for having provided the current best-known values.
Mona Jeunehomme, Sol\`ene Kojtych, Jeffrey Larson, Tom Ragonneau, Viviane Rochon~Montplaisir, Ludovic Salomon, Bastien Talgorn, and students from the MTH8418 derivative-free optimization course at Polytechnique Montr\'eal.

\bibliographystyle{plain}
\bibliography{bibliography}

\pdfbookmark[1]{References}{sec-refs}

\end{document}